\documentclass[11pt]{article}

\usepackage[margin=1in]{geometry}
\usepackage[T1]{fontenc}
\usepackage{amsmath}
\usepackage{amsthm}
\usepackage{mathtools}
\usepackage{amsfonts}
\usepackage{newtxtext}
\usepackage{newtxmath}

\usepackage{algorithm}
\usepackage{algpseudocode}
\usepackage{comment}
\usepackage{subcaption}
\usepackage{multirow}
\usepackage{booktabs}
\usepackage{setspace}
\usepackage{graphicx}
\usepackage{nicefrac}
\usepackage{microtype}
\usepackage{xcolor}
\usepackage[normalem]{ulem}
\usepackage{accents}
\usepackage{chngcntr}

\usepackage{natbib}
\bibpunct[, ]{(}{)}{,}{a}{}{,}%
\usepackage[colorlinks,linkcolor=blue,citecolor=blue,hypertexnames=false]{hyperref}
\colorlet{blue}{black}
\usepackage[nameinlink,capitalize,noabbrev]{cleveref}

\theoremstyle{plain}
\newtheorem{thm}{Theorem}[section]
\newtheorem{lem}{Lemma}[section]
\newtheorem{prop}{Proposition}[section]
\newtheorem{cor}{Corollary}[section]

\theoremstyle{definition}
\newtheorem{assumption}{Assumption}[section]
\newtheorem{defi}{Definition}[section]

\theoremstyle{remark}
\newtheorem{rem}{Remark}[section]

\renewcommand{\proof}[1]{\par\noindent\textit{#1}\ }
\renewcommand{\endproof}{\par\medskip}
\newcommand{\Halmos}{\hfill$\square$}

\DeclareMathOperator*{\argmin}{argmin}

\DeclareMathOperator*{\st}{s.t.}

\newcommand{\rom}[1]{\uppercase\expandafter{\romannumeral #1\relax}}

\newcommand{\CVaR}{\ensuremath{\mathrm{CVaR}}}
\newcommand{\ECSwitch}{\clearpage\appendix}
\newcommand{\ECHead}[1]{\section*{#1}\addcontentsline{toc}{section}{#1}}

\numberwithin{equation}{section}
\begin{document}

\title{Projection-Free Functional Constrained Optimization for Risk Aversion and Sparsity Control}
\author{%
Yi Cheng\thanks{H. Milton Stewart School of Industrial and Systems Engineering, Georgia Institute of Technology, Atlanta, Georgia 30332, USA. Email: \texttt{cheng.yi@gatech.edu}.}
\and Guanghui Lan\thanks{H. Milton Stewart School of Industrial and Systems Engineering, Georgia Institute of Technology, Atlanta, Georgia 30332, USA. Email: \texttt{george.lan@isye.gatech.edu}.}
\and Saeed Masiha\thanks{College of Management of Technology, EPFL, Lausanne, Switzerland. Email: \texttt{mohammadsaeed.masiha@epfl.ch}.}
\and H. Edwin Romeijn\thanks{H. Milton Stewart School of Industrial and Systems Engineering, Georgia Institute of Technology, Atlanta, Georgia 30332, USA. Email: \texttt{edwin.romeijn@isye.gatech.edu}.}
}
\date{}

\maketitle

\begin{abstract}
We study projection-free methods for functional constrained optimization with convex or smooth nonconvex objectives. Such problems arise in applications such as portfolio optimization and radiation therapy planning, where risk-aware criteria and sparsity frequently appear together. For the convex setting, we propose a Level Conditional Gradient (LCG) method that combines a level-set outer loop with a conditional gradient oracle for saddle-point subproblems, and we show an iteration complexity of $\mathcal{O}\big(\epsilon^{-2}\log(\epsilon^{-1})\big)$ for smooth and nonsmooth cases without dependence on the magnitude of an optimal dual Lagrange multiplier. For the nonconvex setting, we propose the Inexact Proximal Point LCG (IPP-LCG) method, which solves a sequence of convex subproblems by LCG and attains $\mathcal{O}\big(\epsilon^{-3}\log(\epsilon^{-1})\big)$ complexity for computing an \((\epsilon,\epsilon)\)-near-KKT point. Numerical results on portfolio selection and IMRT illustrate the practical sparsity/risk trade-offs of the proposed methods.
\end{abstract}


\section{Introduction}\label{sec:intro}
We study functional constrained optimization problems of the form
\begin{equation}
	\begin{aligned}
	f^* := \min\ & f(x)\\
	           \text{s.t.}\ &h_i(x) \le 0,\ i=1,\cdots,m,\\
	          & x\in X,
	\end{aligned}
	\label{cvxctr-model}
\end{equation}
where \(X \subseteq \mathbb{R}^n\) is a nonempty compact convex set, each \(h_i: X \to \mathbb{R}\) is proper, lower semicontinuous, and convex, and \(f: X \to \mathbb{R}\) is proper and lower semicontinuous. We study two regimes: a \emph{convex} setting in which \(f\) is convex (possibly nonsmooth), and a \emph{nonconvex} setting in which \(f\) is differentiable with Lipschitz continuous gradients.

Many motivating instances combine nonlinear aggregate criteria with sparse or limited-actuation decisions. Such models fit naturally into \eqref{cvxctr-model}, but they are challenging to solve at scale by projection-based first-order methods and are also not directly amenable to existing projection-free approaches. In Section~\ref{sec:numerical}, we develop the portfolio and IMRT instances studied in this paper; here we focus on the projection-free methodology for solving \eqref{cvxctr-model}.

Projection-based schemes, including penalty, augmented Lagrangian, and primal--dual methods \citep{lan2013iteration,lan2016iteration,NJLS09,boob2019stochastic,boob2022level}, can handle general functional constraints but require repeated projections onto \(X\), which may be expensive and may destroy sparse structure. Existing projection-free approaches either focus on simpler settings \citep{lan2016conditional,lan2017conditional,qu2018non} or, for general nonlinear constraints, have complexity bounds that depend on the size of an optimal dual multiplier \citep{lan2020conditional}. Our goal is to retain projection-free updates while avoiding such multiplier dependence.

Conditional gradient methods replace projections by calls to a linear minimization oracle (LMO) over \(X\). On polytopes, LMO solutions can be chosen as extreme points, so iterates admit sparse convex-hull representations when the method starts from a sparse point \citep{clarkson2010coresets,jaggi2013revisiting}. This yields an explicit \emph{accuracy--sparsity trade-off}: fewer iterations give sparser representations, while higher accuracy typically requires more iterations.\footnote{For the simplex \(X=\{x\ge 0:\ \sum_i x_i=1\}\), atomic sparsity coincides with coordinate sparsity.} Our aim is to extend projection-free methods to functional constrained problems while retaining the LMO-based iterate structure of conditional-gradient methods, so that the usual finite-iteration accuracy--structure trade-off remains available in practice.

The main difficulty is how to measure progress toward feasibility and optimality in a projection-free way when the constraints are nonlinear. This is substantially harder than in standard Frank--Wolfe settings with simple explicit constraints. The closest related projection-free method is due to \citet{lan2020conditional}, which also handles general nonlinear constraints, but its complexity bounds depend on the magnitude of an optimal dual multiplier. Our approach instead solves compact saddle subproblems that return explicit lower and upper certificates of progress, which drive the level updates and stopping rules while keeping the method projection-free. For the nonconvex setting, we combine this oracle with a proximal-point reduction, which allows us to obtain near-KKT guarantees for the original nonconvex problem while keeping every inner solve projection-free and LMO-based.


	\subsection{Contributions and Outline}
	\paragraph{Contributions.}
	Our main contributions are as follows.
\begin{itemize}
		\item \textbf{LCG for convex functional constraints.} We propose the Level Conditional Gradient (LCG) method for convex instances of \eqref{cvxctr-model}. LCG combines a level-set outer loop with a conditional-gradient oracle (CGO) for saddle subproblems. The oracle returns a primal iterate and a certificate interval \([L_k,U_k]\) that drives the level updates and stopping rule. We show that LCG achieves \(\mathcal{O}(\epsilon^{-2}\log(\epsilon^{-1}))\) CGO-iteration complexity for obtaining an \(\epsilon\)-optimal and \(\epsilon\)-feasible solution, without dependence on an optimal Lagrange multiplier magnitude.

			\item \textbf{IPP-LCG for smooth nonconvex objectives.} For the nonconvex setting, we develop the Inexact Proximal Point LCG (IPP-LCG) method, which solves convex proximal subproblems by LCG. We show that IPP-LCG computes an \((\epsilon,\epsilon)\)-near-KKT point with iteration complexity \(\mathcal{O}(\epsilon^{-3}\log(\epsilon^{-1}))\), while keeping every inner solve projection-free and LMO-based.

			\item \textbf{Empirical evaluations.} We study convex and smooth nonconvex formulations for portfolio selection and IMRT. In portfolio selection, the experiments show that under a common wall-clock budget, LCG and IPP-LCG deliver the strongest feasible out-of-sample step-risk on most datasets while respecting the target support level, relative to CoexDurCG, projection-based first-order baselines, and sparsity-oriented heuristics. In IMRT, the experiments show that on the convex formulation, LCG improves constraint satisfaction relative to CoexDurCG at similar runtimes, and that on the large prostate instance, a sparse convex LCG plan provides an effective warm start for nonconvex IPP-LCG refinement, reducing the maximum listed violation substantially while preserving a clinically sparse treatment structure.
\end{itemize}


\subsection{Related Work}
The functional-constrained template \eqref{cvxctr-model} is classical; the gap addressed here is specifically on the projection-free side. For convex functional constraints, a broad first-order literature uses penalty, augmented Lagrangian, primal--dual, extragradient, and saddle-point schemes \citep{hestenes1969multiplier,powell1969method,lan2013iteration,lan2016iteration,korpelevich76,nemirovski2004prox,NJLS09,boob2019stochastic,hamedani2021primal,boob2022level}. These methods handle general nonlinear constraints and often enjoy strong complexity guarantees, but their iteration models are built around projections or proximal evaluations over \(X\). Our contribution differs from this line at the oracle level: we use only linear minimization over \(X\) and never project onto the full functionally constrained region.

Another ingredient comes from value-function, level-set, and root-finding methodology \citep{lemarechal1995new,nesterov2018lectures,van2009probing,van2011sparse,aravkin2013variational,lin2018level,aravkin2019level}. These works show how to solve optimization problems through a sequence of feasibility or value-evaluation tasks. LCG borrows this outer viewpoint, but the inner solver is different: the feasibility subproblems are handled by a projection-free compact saddle oracle that returns an explicit certificate interval \([L_k,U_k]\), which then drives the level updates and stopping rule.

Projection-free methods themselves go back to the conditional-gradient procedures of \citet{frank1956algorithm} and \citet{levitin1966constrained}; modern analyses emphasizing sparse atomic representations and LMO-based updates include \citet{clarkson2010coresets,jaggi2013revisiting}. Extensions to saddle and constrained settings include \citet{lan2016conditional,lan2017conditional,gidel2017frank,chen2020efficient}. The closest prior work to our convex setting is \citet{lan2020conditional}, which develops CoexCG/CoexDurCG for general affine and nonlinear constraints. Relative to that line, our convex contribution combines a level-set outer loop with an inner compact saddle oracle that returns computable lower and upper certificates, and the resulting complexity bound avoids dependence on the magnitude of an optimal Lagrange multiplier.

For nonconvex problems, one important line uses proximal-point, majorization, penalty, or related subproblem-based ideas to reduce the original problem to a sequence of easier updates \citep{jiang2019structured,guler1992new,scutari2016parallel,lan2019accelerated,kong2019complexity,ma2019proximally,boob2020feasible}. These methods are closest to IPP-LCG in the type of guarantee they seek, but they typically rely on projection/proximal subproblem solvers or on simpler feasible-set models than the general nonlinear functional constraints considered here. A different line studies direct projection-free methods for nonconvex objectives over simple compact domains \citep{qu2018non,oliveira2023note} and, more recently, one-sided or related projection-free methods for nonconvex-concave saddle or CCCP-type reformulations \citep{boroun2023projection,yurtsever2022cccp}. These results are relevant, but they target different models and residuals: simple-set nonconvex optimization or stationarity/saddle guarantees for reformulated minimax problems. Our nonconvex result instead keeps the original nonlinear constraints, solves each convexified constrained subproblem with the projection-free LCG routine, and establishes an \((\epsilon,\epsilon)\)-near-KKT guarantee for the original functionally constrained problem.

\subsection{Notations}
The following notations will be used throughout the paper. Without specific mention, $\|\cdot\|$ denotes an arbitrary norm (not necessarily associated with an inner product) and $\|\cdot\|_{*}$ denotes its conjugate. The notation $[1:n]$ represents the set of integers $\{1, \dots, n\}$. For a closed convex set $X\subset \mathbb{R}^n$, we write \(\operatorname{ri}(X)\) for its relative interior. The set $N_X(x)$ denotes the normal cone at $x\in X$ and $N_X(x):=\{g\in \mathbb{R}^n| \forall z\in X: g^\top(z-x)\le 0\}$. A function $f: \mathbb{R}^n \rightarrow \mathbb{R}$ is $L_f$-smooth if $\|\nabla f(x_1) - \nabla f(x_2)\|_* \le L_f\|x_1 - x_2\|$, $\forall x_1,x_2\in X$. A function $f: \mathbb{R}^n \rightarrow \mathbb{R}$ is $M_f$-Lipschitz continuous if $|f(x_1) - f(x_2)|\le M_f\|x_1 - x_2\|$, $\forall x_1,x_2\in X$. Let \( x^* \) be an optimal solution of \eqref{cvxctr-model}, and define \( f^* := f(x^*) \). A point \( \bar{x} \) is an \(\epsilon\)-optimal and \(\epsilon\)-feasible solution (or \(\epsilon\)-solution) of \eqref{cvxctr-model} if \( \bar{x} \in X \), \( f(\bar{x}) - f^* \leq \epsilon \), and \(\max\limits_{i \in [1:m]} h_i(\bar{x}) \leq \epsilon\). We use $x_{t,i}$ to denote the $i$-th element of vector $x_t$.

Section~\ref{sec-level-cg} presents LCG and CGO for convex problems, Section~\ref{sec-nonconvex} introduces IPP-LCG for nonconvex problems, and Section~\ref{sec:numerical} combines the applications and numerical results. Additional proofs and supporting experiments are given in the e-companion.

\section{Convex Conditional Gradient Method}\label{sec-level-cg}
In this section, we introduce the Level Conditional Gradient (LCG) method that solves the convex constrained optimization problem \eqref{cvxctr-model} and establish its convergence rate.
We start by stating the standing assumptions on the convex constrained model.
\begin{assumption}\label{asm:lcg-problem}
For problem \eqref{cvxctr-model}, \(f\) is proper, lower semicontinuous, convex, and \(M_f\)-Lipschitz continuous on \(X\). Each \(h_i\) is proper, lower semicontinuous, convex, and \(M_{h_i}\)-Lipschitz continuous on \(X\).
\end{assumption}

It is well-known that problem \eqref{cvxctr-model} can be reduced to a root finding problem. 
For a given level estimate $l \in \mathbb{R}$, let us define 
 \begin{align}
\phi(l) &:= \min\limits_{x\in X} \max\left\{ f(x) -l, h_1(x), \ldots, h_m(x) \right\}\nonumber\\
&=\min_{x \in X} \max\limits_{(\gamma,z) \in Z}\ \gamma[ f(x) - l] + \sum\limits_{i=1}^m z_i h_i(x). \label{cvxctr-model3}
\end{align}
Here
$Z := \{(\gamma, z) \in \mathbb{R}^{m+1}: \gamma + \sum_{i=1}^m z_i = 1, \gamma, z_i\ge 0\}$ denotes the standard simplex. 
We can easily verify that: (a) $\phi(l)$ is monotonically non-increasing and convex w.r.t. $l$; (b) $\phi(f^*) =0$; (c) $\phi(l) \ge 0$ for any $l \le f^*$ and $\phi(l) \le 0$ for any $l \ge f^*$. Therefore, problem \eqref{cvxctr-model} is equivalent to
finding the root of $\phi(l) = 0$.

We propose to solve \eqref{cvxctr-model} by LCG (see Algorithm \ref{alg:level-CG}), which consists of an \textit{outer loop} that updates the
level estimate $l$ (i.e., the estimation of $f^*$),
and an \textit{inner loop} that calls a specialized conditional gradient oracle (CGO) to solve the saddle point problem in \eqref{cvxctr-model3} given a level estimate $l$.

To guide the discussion, we first describe the outer level-update mechanism. We then present the inner CGO routine for smooth saddle subproblems and derive the corresponding overall complexity of LCG. After that, we extend the same construction to the structured nonsmooth case and state its analogous overall complexity bound.
\subsection{Outer Loop of LCG}
The basic idea of the LCG method is to apply an approximate Newton's method to solve $\phi(l) =0$. Assume for the moment that problem \eqref{cvxctr-model3} can be solved exactly for a given $l$ ($l=l_k$). Then one can compute the function value $\phi(l_k)$ and a subgradient $\phi'(l_k)$. Solving the linear equation $\phi(l_k) + \phi'(l_k) (l - l_k) = 0$ gives the updated iterate $l_{k+1} = l_k - \tfrac{\phi(l_k)}{\phi'(l_k)}$.
Since $\phi(l_k)$ cannot be computed exactly, we suggest to use a computable lower bound and an approximate subgradient in place of $\phi(l_k)$ and $\phi'(l_k)$ in the above equation, respectively. Started with an initial level estimate $l_1\le f^*$, we call CGO to compute a lower bound $L_k$, an upper bound $U_k$ of $\phi(l_k)$
and an approximate primal-dual pair $(x_{k}, \bar z_k) \in X \times Z$ of problem \eqref{cvxctr-model3} at the $k$-th iteration (see Algorithm \ref{alg:level-CG}). For this specialized simplex dual set, we write $\bar z_k = (\gamma_k, \zeta_k)$, so $\gamma_k$ is the first dual coordinate associated with the level term $f(x)-l_k$.
A gap defined by these bounds (i.e. $U_k - L_k$) indicates how accurately problem \eqref{cvxctr-model3} (with $l=l_k$) is solved.
Whenever the upper bound $U_k \le \epsilon$, LCG terminates since an approximate root of $\phi(l) = 0$ has been found due to $\phi(l_k) \le U_k \le \epsilon$ and $\phi(l_k) \ge 0$.
Otherwise, the algorithm updates the level estimate $l_k$. 
More specifically, we define the following linear function as a lower approximation of $\phi(l),\forall l\in \mathbb{R}$:
\begin{equation}
\mathcal{L}_k(l):=L_k -\gamma_k(l-l_k).
    \label{lower-bounding-L}
\end{equation}
 Intuitively, $\mathcal{L}_k(l)$ underestimates $\phi(l)$ since $-\gamma_k$ and $L_k$ respectively serve as an approximate subgradient and a lower bound for $\phi(\cdot)$ at $l=l_k$. 
To perform the approximate Newton's step as mentioned earlier, we solve $\mathcal{L}_k(l) = 0$ and obtain the following update of the level estimate 
 \begin{equation} \label{eq:updating_level}
 	l_{k+1} = l_k + \frac{1}{\gamma_k}L_k.
 \end{equation}
We note that the LCG method provides a general framework for solving the root finding problem in \eqref{cvxctr-model3} and is not restricted to a particular inner oracle, as long as the output 	
($\gamma_k, L_k,U_k$) of the inner oracle (e.g., CGO) satisfies the following conditions: 
\begin{align}
   \gamma_k&>0, \label{eq:c1}\\
L_k \le &\phi(l_k) \le U_k, \label{eq:c2}\\
\mathcal{L}_k(l) &\le \phi(l), \forall l. \label{eq:c3}
\end{align}
		\begin{algorithm}[h]
			\caption{Level Conditional Gradient Method (LCG)} 
			\begin{algorithmic}[1]
				\State {Inputs: $\epsilon >0, \mu\in(\frac{1}{2},1)$.}
				\State {Initialization: choose $x_0\in X$, pick $g_0\in\partial f(x_0)$, and set $l_1 = \min\{ f(x_0) + \langle g_0,x-x_0\rangle:\ {x\in X} \}$.}\label{alg: init}
				\For {$k=1,2,\ldots$}
				\State {Call CGO with input $l_k$ and obtain $(x_k,\bar z_k)$, lower bound $L_k$, upper bound $U_k$ such that $U_k- L_k\le (1-\mu)\epsilon$.}
				\State {Set $\gamma_k = (\bar z_k)_0$.}
				\If {$U_k\le \epsilon$}
			\State{Terminate and return $x_{k}$.}
			\EndIf
			\State {$l_{k+1} = l_k + \frac{1}{\gamma_{k}} L_k$.} \label{alg: level-CG-update t}
			\EndFor
	\end{algorithmic} 
	\label{alg:level-CG}
\end{algorithm}
The following lemma states an important property of the sequence of the level estimates $(l_k)_{k\ge 1}$ generated in the outer loop of the algorithm. 
Such property will be used in establishing the number of outer loops required by the LCG method.
		\begin{lem}
		At  iteration $k$, if Algorithm \ref{alg:level-CG} does not terminate, then $L_k > 0$. Moreover, the sequence of the level estimates satisfies $l_1< \cdots < l_k < l_{k+1}< \cdots\le f^*$, $k\ge 1$. Consequently, \(0=\phi(f^*) \le \phi(l_{k+1}) \le \phi(l_k)\le \cdots \le \phi(l_1)\).
		\label{lem5}
		\end{lem}
	\proof{Proof.}
	If Algorithm \ref{alg:level-CG} does not terminate at iteration \(k\), then \(U_k>\epsilon\). On the other hand, the stopping rule of CGO guarantees \(U_k-L_k\le (1-\mu)\epsilon\). Hence \(L_k \ge U_k-(1-\mu)\epsilon > \epsilon-(1-\mu)\epsilon = \mu\epsilon >0\). Since \(\gamma_k>0\) by \eqref{eq:c1}, the level update \eqref{eq:updating_level} gives \(l_{k+1}-l_k=\tfrac{L_k}{\gamma_k}>0\), so the level sequence is strictly increasing.

	Next, by construction of \(l_{k+1}\), we have \(\mathcal L_k(l_{k+1})=0\). Because \(\mathcal L_k(\cdot)\) underestimates \(\phi(\cdot)\) by \eqref{eq:c3}, it follows that \(0=\mathcal L_k(l_{k+1})\le \phi(l_{k+1})\). At the same time, \(\phi\) is nonincreasing and \(\phi(f^*)=0\). Therefore \(\phi(l_{k+1})\ge 0=\phi(f^*)\) implies \(l_{k+1}\le f^*\). Since \(l_1\le f^*\) by initialization, we obtain \(l_1<l_2<\cdots<l_k<l_{k+1}\le f^*\). Finally, because \(\phi\) is nonincreasing and \(l_k\le l_{k+1}\le f^*\), we have \(0=\phi(f^*)\le \phi(l_{k+1})\le \phi(l_k)\le \cdots \le \phi(l_1)\).
	This proves the claim. \Halmos
	\endproof
\begin{rem}
In LCG, we require that the initial level estimate satisfies $l_1 \le f^*$. Otherwise, the algorithm terminates at the first outer iteration. To see this, if $l_1 \ge f^*$, then $L_1 \le \phi(l_1) \le \phi(f^*) =0$, which holds because $L_1$ is the lower bound of $\phi(l_1)$ and $\phi$ is a non-increasing function. Hence by the stopping criteria of CGO, $U_1 \le L_1 + (1-\mu)\epsilon \le \epsilon$, which results in the termination of the algorithm. A decreasing-level variant can also approximate $f^*$ from above when started with $l_1 \ge f^*$, but we do not pursue that variant here.
\label{rem-outer-opposite-update}
\end{rem}
In the theorem below, we establish the iteration complexity of reaching ``$U_k \le \epsilon$'', which is essentially the outer loop iteration complexity of solving \eqref{cvxctr-model} by Algorithm \ref{alg:level-CG}. 
		\begin{thm}\label{prop2}
				For all $k\ge 1$, we have
				\begin{equation}
				U_k \le \left( f^* - l_1\right)\frac{1}{\mu}\left(\frac{1}{2\mu}\right)^{k-1},
			\label{prop-outer-complexity}
				\end{equation} where $\mu \in (\frac{1}{2}, 1)$, $l_1$ is the initial estimate of the optimal value of (\ref{cvxctr-model}) such that $l_1 \le f^*$.
				Moreover, given precision $\epsilon$, at the termination of LCG when $U_k\le \epsilon$, the algorithm yields an $\epsilon$-optimal and $\epsilon$-feasible solution $x_k$ of problem \eqref{cvxctr-model}. 
			\end{thm}
	\proof{Proof.}
	Fix an iteration \(k\) at which the algorithm has not yet terminated. Since \(\mathcal L_k(\cdot)\) is affine, \(\mathcal L_k(l_k)=L_k\), and \(\mathcal L_k(l_{k+1})=0\), equality of slopes gives \([\mathcal L_k(l_{k-1})-L_k]/(l_k-l_{k-1})=L_k/(l_{k+1}-l_k)\). Hence \(\mathcal L_k(l_{k-1})=L_k[1+(l_k-l_{k-1})/(l_{k+1}-l_k)]\). Multiplying by \(l_{k+1}-l_k\) and using \(a+b\ge 2\sqrt{ab}\) for \(a,b\ge 0\), we obtain \((l_{k+1}-l_k)\mathcal L_k(l_{k-1})=L_k[(l_{k+1}-l_k)+(l_k-l_{k-1})]\ge 2L_k\sqrt{(l_{k+1}-l_k)(l_k-l_{k-1})}\).
	Therefore,
	\begin{equation}
	\frac{\mathcal L_k(l_{k-1})}{\sqrt{l_k-l_{k-1}}}
	\ge
	\frac{2L_k}{\sqrt{l_{k+1}-l_k}}.
	\label{eq:outer-proof-key}
	\end{equation}

	Since the algorithm has not terminated at iteration \(k\), we have \(U_k>\epsilon\) and \(U_k-L_k\le (1-\mu)\epsilon\). Because \(\epsilon<U_k\), this yields \(U_k-L_k \le (1-\mu)U_k\), hence \(L_k\ge \mu U_k\). Also, \(U_{k-1}\ge \phi(l_{k-1})\ge \mathcal L_k(l_{k-1})\) by \eqref{eq:c2}--\eqref{eq:c3}. Combining these inequalities with \eqref{eq:outer-proof-key} gives
	\(U_{k-1}\ge \mathcal L_k(l_{k-1})\ge 2\mu U_k\sqrt{l_k-l_{k-1}}/\sqrt{l_{k+1}-l_k}\),
	or equivalently
	\begin{equation}
	U_k \le \frac{1}{2\mu}\sqrt{\frac{l_{k+1}-l_k}{l_k-l_{k-1}}}\,U_{k-1}.
	\label{eq:outer-recursion}
	\end{equation}

	Applying \eqref{eq:outer-recursion} recursively yields \(U_k \le \left(\tfrac{1}{2\mu}\right)^{k-1}\sqrt{\tfrac{l_{k+1}-l_k}{l_2-l_1}}\,U_1\). Next, \(U_1-L_1\le (1-\mu)\epsilon\) implies \(U_1\le L_1/\mu\) whenever the algorithm reaches a second iteration; if it terminates at the first iteration, the theorem is immediate. Since \(L_1\le \phi(l_1)\) and \(\phi\) is \(1\)-Lipschitz with respect to the level argument, we have \(\phi(l_1)\le f^*-l_1\). Also \(\mathcal L_1(l_2)=0\) and \(\mathcal L_1(l_1)=L_1\), so \(l_2-l_1=L_1/\gamma_1\ge L_1\) because \(\gamma_1\le 1\). Finally, Lemma \ref{lem5} gives \(l_{k+1}-l_k\le f^*-l_1\). Substituting these bounds into the previous display gives \(U_k \le \tfrac{1}{\mu}\left(\tfrac{1}{2\mu}\right)^{k-1}(f^*-l_1)\), which is \eqref{prop-outer-complexity}. For the final claim, if the algorithm terminates at iteration \(k\), then \(U_k\le \epsilon\). Since \(\phi(l_k)=\min_{x\in X}\max\{f(x)-l_k,h_1(x),\ldots,h_m(x)\}\le U_k\), the returned point \(x_k\) satisfies \(f(x_k)-l_k\le U_k\) and \(\max_{i\in[m]} h_i(x_k)\le U_k\). Because \(l_k\le f^*\) by Lemma \ref{lem5}, we obtain \(f(x_k)-f^*\le f(x_k)-l_k\le U_k\le \epsilon\) and \(\max_{i\in[m]} h_i(x_k)\le \epsilon\). Hence \(x_k\) is \(\epsilon\)-optimal and \(\epsilon\)-feasible. \Halmos
	\endproof


\subsection{Conditional Gradient Oracle for Smooth Saddle Subproblems}
\label{sec-smooth-functions}

We first present CGO for smooth saddle problems, since this is the clearest setting in which to see the oracle structure and the certificate mechanism. The structured nonsmooth case is then obtained by smoothing and minor modifications described in Section~\ref{sec-cgo-nonsmooth}. The variables \(x_t,p_t\) play the primal roles, \(r_t,\bar z_t\) are dual iterates, and \(L_t,U_t\) are computable lower/upper certificates for the saddle subproblem.
CGO solves the saddle-point problem
\begin{align}
 \begin{split}
\bar \phi := \min\limits_{x\in \bar{X}} \max\limits_{\bar z\in \bar{Z}}\ \bar{f}(x) + \sum\limits_{i=1}^{\bar m} \bar z_i \bar{h}_i(x).
 \end{split}
 \label{model}
\end{align}
Here, $\bar{f}: \bar{X} \rightarrow \mathbb{R}$ and $\bar{h} : \bar{X} \rightarrow \mathbb{R}^{\bar{m}}$ are proper lower semicontinuous convex functions, $ \bar{X}\subset \mathbb{R}^{\bar{n}}$ is a nonempty compact convex set, and $\bar{Z} \subset \mathbb{R}_{+}^{\bar{m}}$ is a convex compact set. 
Under these assumptions an optimal pair of solutions $(x^*,\bar z^*) \in \bar{X}\times \bar{Z}$ of problem (\ref{model}) must exist.
Clearly, 
the subproblem in (\ref{cvxctr-model3}) can be viewed as 
a special case of problem \eqref{model} with $\bar{f} = 0$, $\bar{X} = X$, $\bar{Z} = Z$ and $\bar{h} = (f-l, h)$.
The smooth analysis below relies on the following regularity conditions.
\begin{assumption}\label{asm:cgo-smooth}
For problem \eqref{model}, \(\bar{f}\) is convex, \(L_{\bar f}\)-smooth, and \(M_{\bar f}\)-Lipschitz continuous over \(\bar X\). For each \(i=1,\ldots,\bar m\), \(\bar h_i\) is convex, \(L_{\bar h_i}\)-smooth, and \(M_{\bar h_i}\)-Lipschitz continuous over \(\bar X\).
\end{assumption}
Let $\nu: \bar{Z}\rightarrow \mathbb{R}$ be a $1$-strongly convex and $L_{\nu}$-smooth distance generating function and define the proximal function at point $z'\in \bar{Z}$ as $V(z',z) := \nu(z) - \nu(z') - \langle \nabla \nu(z'), z - z'\rangle$ for $z\in \bar{Z}$. Further, denote the linear approximation of $\bar{f}$ and $\bar{h}$ at $x'$ as $\ell_{\bar{f}}(x',x) := \bar{f}(x') + \langle \nabla\bar{f}(x'), x - x' \rangle$ and $\ell_{\bar{h}_i}(x',x) :=  \bar{h}_i(x') + \langle \nabla\bar{h}_i(x'), x - x' \rangle$ for \(i =1,\cdots,\bar{m}\).
The algorithmic scheme of CGO is stated in Algorithm \ref{alg:CG oracle}. It updates dual variables through extrapolation and proximal steps (\eqref{alg-cg-extrap}--\eqref{alg-cg-z}), computes the primal direction by one LMO call (\eqref{alg-cg-primal}), and maintains lower/upper certificates \(L_t,U_t\) through \eqref{alg-cg-lower}--\eqref{alg-cg-upper}.  

		\begin{algorithm}[h]
			\caption{Conditional Gradient Oracle (CGO) } 
			\begin{algorithmic}[1]
				\State{Parameters: $\lambda_t\ge 0, \tau_t \ge 0, \alpha_t\in [0,1], \alpha_1 = 1, \epsilon >0, \mu \in (\frac{1}{2},1)$.}
					\State {Initialization: $x_{-2} = x_{-1} = x_0 \in \bar{X}$, $p_{-1} = p_0 \in \bar{X}$,  $z_0 = r_0\in \bar{Z}$, $\underline{f}_0(x)=\ell_{\bar{f}}(x_0,x)$, $\underline{h}_0(x)=\langle \ell_{\bar{h}}(x_0,x), r_0\rangle$.}
				\For {$t=1,2,\ldots$}
				\State Compute $\bar z_t$, $x_t$, $L_t$ and $U_t$ according to
				\State $\tilde{h}_t = \ell_{\bar{h}}(x_{t-2},p_{t-1}) + \lambda_t[\ell_{\bar{h}}(x_{t-2},p_{t-1}) - \ell_{\bar{h}}(x_{t-3},p_{t-2})]$. \label{alg-cg-extrap}
				\State $r_t = \argmin\limits_{\bar z\in \bar{Z}}\ \langle -\tilde{h}_t, \bar z \rangle + \tau_tV(r_{t-1},\bar z)$. \label{alg-cg-dual}
				\State $\bar z_t = (1-\alpha_t)\bar z_{t-1} + \alpha_tr_t$. \label{alg-cg-z}
				\State $p_t = \argmin\limits_{x\in \bar{X}}\ \ell_{\bar{f}}(x_{t-1},x) + \langle \ell_{\bar{h}}(x_{t-1},x), r_t \rangle$. \label{alg-cg-primal}
				\State $x_t = (1-\alpha_t)x_{t-1} + \alpha_tp_t$. \label{alg-cg-x}
				\State $\underline{f}_t(x) = (1-\alpha_t)\underline{f}_{t-1}(x)  + \alpha_t\ell_{\bar{f}}(x_{t-1},x)$.
					\State $\underline{h}_t(x)  = (1-\alpha_t)\underline{h}_{t-1}(x) + \alpha_t\langle \ell_{\bar{h}}(x_{t-1},x), r_t \rangle$.
					\State $L_t = \min\limits_{x\in \bar{X}}\ \underline{f}_t(x) + \underline{h}_t(x)$. \label{alg-cg-lower}
					\State $U_t = \max\limits_{\bar z\in \bar{Z}}\ \bar{f}(x_t) + \langle \bar{h}(x_t),\bar z \rangle$. \label{alg-cg-upper}
				  \If {$U_t - L_t \le (1-\mu)\epsilon$}
			   \State {Terminate and return $x_t, \bar z_t, L_t, U_t$.}
			   \EndIf
			   \EndFor
			\end{algorithmic} 
			\label{alg:CG oracle}
			\end{algorithm}
\begin{equation}
\underline{f}_t(x) = (1-\alpha_t)\underline{f}_{t-1}(x)  + \alpha_t\ell_{\bar{f}}(x_{t-1},x).
\label{alg-cg-lowerf}
\end{equation}
\begin{equation}
\underline{h}_t(x) = (1-\alpha_t)\underline{h}_{t-1}(x) + \alpha_t\langle \ell_{\bar{h}}(x_{t-1},x), r_t \rangle.
\label{alg-cg-lowerh}
\end{equation}
For notational simplicity, we write \(z_t:=\bar z_t\) in the convergence discussion below and let \(\bar M:=\left(\sum_{i=1}^{\bar m}M_{\bar h_i}^2\right)^{1/2}\), \(D_{\bar X}:=\max_{x_1,x_2\in \bar X}\|x_1-x_2\|\), and \(\bar V:=\max_{z_1,z_2\in \bar Z}V(z_1,z_2)\).
The smooth CGO parameters are chosen as
\begin{equation}
\alpha_t=\frac{2}{t+1},\qquad \lambda_t=\frac{t-1}{t},\qquad \tau_t=9\sqrt{t}\,\bar M D_{\bar X},\qquad t\ge 1.
\label{thm1-params1}
\end{equation}

It is worth briefly comparing CGO with CoexCG/CoexDurCG \citep{lan2020conditional}. The main differences relevant here are the following. First, when CGO is used inside LCG for \eqref{cvxctr-model3}, the dual set is the simplex \(\bar Z=Z\), so the resulting LCG complexity bound does not depend on the magnitude of an optimal Lagrange multiplier. Second, CGO maintains explicit computable certificates \(L_t\) and \(U_t\) and terminates when this certificate gap is small. Third, because \(\bar Z\) is compact in our setting, the smooth-case choice of \(\tau_t\) in \eqref{thm1-params1} does not need prior knowledge of the horizon; see Remark \ref{rmk_depend_total_iterations}. Moreover, the gap function for CGO is based on the upper bound \(U_t\) in \eqref{alg-cg-upper} and the computable lower bound \(L_t=\min_{x\in\bar{X}}\underline f_t(x)+\underline h_t(x)\) in \eqref{alg-cg-lower}. Since this lower certificate is more conservative than the usual benchmark \( \bar f(x)+\langle \bar h(x), z_t\rangle\), controlling the resulting gap yields a stronger certificate than in standard primal--dual analyses with comparable iteration rates.
In the remaining part of this subsection, we discuss the convergence
properties of CGO.\\
The following lemma shows that $L_t$ and $U_t$, $t=1,2, \ldots$ are valid lower bounds and upper bounds of (\ref{model}), respectively. 
\begin{lem}
Let $\bar \phi$, $L_t$ and $U_t$
be defined in \eqref{model}, \eqref{alg-cg-lower} and \eqref{alg-cg-upper}, respectively.
Also let $z_t$ be defined in (\ref{alg-cg-z}). Then, for any $t \ge 1$, we have \(\underline{f}_t(x)\le \bar{f}(x)\), \(\underline{h}_t(x)\le \langle \bar{h}(x), z_t\rangle\), and \(L_t \le \bar \phi \le U_t\).
\label{lem1}
\end{lem}
	\proof{Proof.}
		We first show that \(\underline f_t\) remains a lower model of \(\bar f\). This is true at \(t=0\) because \(\underline f_0(x)=\ell_{\bar f}(x_0,x)\le \bar f(x)\) by convexity of \(\bar f\). If \(\underline f_{t-1}(x)\le \bar f(x)\), then by \eqref{alg-cg-lowerf}, \(\underline f_t(x)=(1-\alpha_t)\underline f_{t-1}(x)+\alpha_t\ell_{\bar f}(x_{t-1},x)\le (1-\alpha_t)\bar f(x)+\alpha_t\bar f(x)=\bar f(x)\). Thus \(\underline f_t(x)\le \bar f(x)\) for all \(t\).

	Next we prove the bound for \(\underline h_t\). Define \(\Gamma_1:=1\) and \(\Gamma_t:=(1-\alpha_t)\Gamma_{t-1}\) for \(t\ge 2\). Dividing \eqref{alg-cg-lowerh} by \(\Gamma_t\) gives
	\[
	\frac{\underline h_t(x)}{\Gamma_t}
	=
	\frac{\underline h_{t-1}(x)}{\Gamma_{t-1}}
	+
	\frac{\alpha_t}{\Gamma_t}\langle \ell_{\bar h}(x_{t-1},x),r_t\rangle.
	\]
		Since each component of \(r_t\) is nonnegative and \(\ell_{\bar h_i}(x_{t-1},x)\le \bar h_i(x)\) by convexity of \(\bar h_i\), we have \(\langle \ell_{\bar h}(x_{t-1},x),r_t\rangle \le \langle \bar h(x),r_t\rangle\).
	Iterating the previous recursion therefore yields
	\[
	\frac{\underline h_t(x)}{\Gamma_t}
	\le
	\sum_{j=1}^t \frac{\alpha_j}{\Gamma_j}\langle \bar h(x),r_j\rangle.
	\]
	Multiplying both sides by \(\Gamma_t\), and defining \(\theta_j^t:=\Gamma_t\alpha_j/\Gamma_j\), we obtain
	\[
	\underline h_t(x)\le \sum_{j=1}^t \theta_j^t\langle \bar h(x),r_j\rangle
	=
	\left\langle \bar h(x),\sum_{j=1}^t\theta_j^t r_j\right\rangle.
	\]
		Because \(z_t=(1-\alpha_t)z_{t-1}+\alpha_t r_t\), the coefficients \(\{\theta_j^t\}_{j=1}^t\) form a convex combination and satisfy \(z_t=\sum_{j=1}^t\theta_j^t r_j\). Hence \(\underline h_t(x)\le \langle \bar h(x),z_t\rangle\).

	Combining the two lower-model bounds gives
	\[
	\underline f_t(x)+\underline h_t(x)\le \bar f(x)+\langle \bar h(x),z_t\rangle
	\le
	\max_{z\in \bar Z}\big\{\bar f(x)+\langle \bar h(x),z\rangle\big\},
	\qquad \forall x\in \bar X.
	\]
		Taking the minimum over \(x\in\bar X\) shows that \(L_t=\min_{x\in\bar X}\{\underline f_t(x)+\underline h_t(x)\}\le \bar \phi\).
	On the other hand, since \(x_t\in\bar X\),
	\[
	U_t
	=
	\max_{z\in \bar Z}\big\{\bar f(x_t)+\langle \bar h(x_t),z\rangle\big\}
	\ge
	\min_{x\in \bar X}\max_{z\in \bar Z}\big\{\bar f(x)+\langle \bar h(x),z\rangle\big\}
	=
	\bar \phi.
	\]
	Therefore \(L_t\le \bar \phi\le U_t\), completing the proof. \Halmos
	\endproof
In view of Lemma~\ref{lem1}, $(\underline{f}_t+ \underline{h}_t)(\cdot)$ provides a lower bound for the objective of \eqref{model}, i.e., $(\underline{f}_t+ \underline{h}_t)(\cdot) \le \bar{f}(\cdot) + \langle \bar{h}(\cdot), z_{t}\rangle$. This motivates us to
define the gap function for problem \eqref{model} as
\begin{equation}
\bar{Q}_t(w_t,w) := \bar{f}(x_t) + \langle \bar{h}(x_t), z\rangle- \underline{f}_t(x) - \underline{h}_t(x), \label{gap}
\end{equation}
where $w_t := (x_t,z_t)$, $w := (x,z)$. 
Also, by the definition of $L_t$ and $U_t$, we can easily see that \(\max\limits_{w\in\bar{X}\times\bar{Z}}\bar{Q}_t(w_t,w) = U_t - L_t\). It is worth mentioning here that the gap function in \eqref{gap} is \emph{different} from the standard primal--dual gap used in the saddle-point literature and in related projection-free work such as \citet{nemirovski2004prox,lan2020conditional}, namely \(\tilde{Q}_t(w_t,w) := \bar{f}(x_t) + \langle \bar{h}(x_t), z\rangle - \bar{f}(x) - \langle \bar{h}(x), z_t\rangle\).
 As a consequence, these algorithms require the solution of $\min_{x \in \bar X} (\bar{f}(x) + \langle \bar{h}(x), z_t\rangle)$ to compute a lower bound on $\bar \phi$, which can be computationally expensive unless both $\bar f$ and $\bar h$ are simple enough (e.g., linear functions). 
 On the other hand, the computation of the lower bound $L_t$ in CGO only requires one call to the linear minimization oracle. In addition, since $\bar{Q}_t(w_t,w) \ge \tilde{Q}_t(w_t,w)$,
 we obtain stronger convergence guarantees for the developed algorithm by using $\bar{Q}_t(w_t,w)$ instead of $\tilde{Q}_t(w_t,w)$
as the error measure.
 
We now state the smooth convergence guarantee for CGO.
\begin{thm}\label{thm:cgo-smooth-main}
Under Assumption \ref{asm:cgo-smooth}, suppose the parameters are set as in \eqref{thm1-params1}.
Then for any \(t\ge 1\) and \(w=(x,z)\in \bar X\times \bar Z\),
\[
\bar Q_t(w_t,w)\le \frac{2(L_{\bar f}+z^\top L_{\bar h})D_{\bar X}^2}{t+1}
+ \frac{\bar M D_{\bar X}}{\sqrt{t+1}}\left(18\bar V+\frac76\right),
\]
where \(L_{\bar h}:=(L_{\bar h_1},\ldots,L_{\bar h_{\bar m}})\). Consequently, to attain \(U_t-L_t\le (1-\mu)\epsilon\), CGO needs \(\mathcal{O}(\epsilon^{-2})\) iterations.
	\end{thm}
	\proof{Proof.}
		Although the certificate \(\bar Q_t\) differs from the standard primal--dual gap, its updates satisfy the same accelerated primal--dual recursion structure as in \citet[Section~3]{lan2020first}; this is stated precisely in Lemma \ref{lem:smooth-gap-recursion} in Appendix~\ref{sec:appen_smooth}. Using this recursion, for every \(T\ge 1\) and \(w=(x,z)\in \bar X\times \bar Z\),
		\begin{equation}
		\bar{Q}_T(w_T,w) \le \Gamma_T\sum_{t=1}^T\left[ \frac{(L_{\bar{f}} + z^{\top}L_{\bar{h}})\alpha_t^2}{2\Gamma_t}D_{\bar{X}}^2 + \frac{9\alpha_t\lambda_t^2\bar{M}^2D_{\bar{X}}^2}{2\tau_t\Gamma_t}\right] + \frac{9\alpha_T\bar{M}^2D_{\bar{X}}^2}{2\tau_T} +\alpha_T\tau_T\bar{V}.
		\label{eq:smooth-gap-recursion-main}
		\end{equation}
		The main change from \citet[Section~3]{lan2020first} is that the terminal dual term is uniformly bounded because \(\bar Z\) is compact.
		Now substitute the parameter choice \(\alpha_t=2/(t+1)\), \(\lambda_t=(t-1)/t\), and \(\tau_t=9\sqrt{t}\bar M D_{\bar X}\). For this choice, $\Gamma_T=2/(T(T+1))$ and ${\alpha_t}/{\Gamma_t}=t$. Hence
	\[
	\Gamma_T\sum_{t=1}^T \frac{\alpha_t^2}{2\Gamma_t}
	=
	\Gamma_T\sum_{t=1}^T \frac{t}{2}\alpha_t
	=
	\Gamma_T\sum_{t=1}^T \frac{t}{t+1}
	\le \frac{2}{T+1}.
	\]
	Also,
	\[
	\Gamma_T\sum_{t=1}^T \frac{9\alpha_t\lambda_t^2}{2\tau_t\Gamma_t}
	=
	\Gamma_T\sum_{t=1}^T \frac{(t-1)^2}{2t^{3/2}\bar M D_{\bar X}}
	\le
	\Gamma_T\sum_{t=1}^T \frac{\sqrt t}{2\bar M D_{\bar X}}
	\le
	\frac{2\sqrt T}{3(T+1)\bar M D_{\bar X}},
	\]
	where we used \(\sum_{t=1}^T \sqrt t \le \frac{2}{3}T^{3/2}\). Finally,
	\[
	\frac{9\alpha_T\bar M^2D_{\bar X}^2}{2\tau_T}
	=
	\frac{\bar M D_{\bar X}}{(T+1)\sqrt T}
	\le
	\frac{\bar M D_{\bar X}}{\sqrt{T+1}},
	\]
	and
	\[
	\alpha_T\tau_T\bar V
	=
	\frac{18\sqrt T}{T+1}\bar M D_{\bar X}\bar V
	\le
	\frac{18\bar M D_{\bar X}\bar V}{\sqrt{T+1}}.
	\]
	Substituting these estimates into \eqref{eq:smooth-gap-recursion-main} yields
	\[
	\bar{Q}_T(w_T,w) \le \frac{2(L_{\bar f}+z^\top L_{\bar h})D_{\bar X}^2}{T+1}
	+ \frac{\bar M D_{\bar X}}{\sqrt{T+1}}\left(18\bar V+\frac76\right),
	\]
	which is the desired bound after renaming \(T\) as \(t\). Since \(U_t-L_t=\max_{w\in\bar X\times \bar Z}\bar Q_t(w_t,w)\), the final complexity statement follows immediately. \Halmos
\endproof

For the LCG subproblem \eqref{cvxctr-model3}, Theorem \ref{thm:cgo-smooth-main} yields
\begin{equation}
U_t-L_t \le \frac{2\max\{L_f,\max_{i\in[m]}L_{h_i}\}D_X^2}{t+1}
+ \frac{\left(M_f^2+\sum_{i=1}^m M_{h_i}^2\right)^{1/2}D_X}{\sqrt{t+1}}
\left[18\bar V+\frac{7}{6}\right].
\label{eq:lcg-smooth-inner-bound}
\end{equation}
Hence each smooth inner CGO call still needs \(\mathcal{O}(\epsilon^{-2})\) iterations. We record the corresponding explicit specialization bounds on \(Q_t(w_t,w)\), \(f(x_t)-f^*\), and \(\max_{i\in[m]} h_i(x_t)\) in Corollary \ref{cor1} in Appendix~\ref{sec:append-lcg}.

\subsection{Overall Complexity in the Smooth Case}\label{sec:level-cg-overall}
We now combine the smooth CGO rate with the outer LCG update. For the subproblem \eqref{cvxctr-model3}, we have \(\bar f(x)=0\), \(\bar h(x;l)=\left(f(x)-l,h(x)\right)\), \(\bar X=X\), and \(\bar Z=Z\). In this specialization, each dual iterate returned by CGO is written as \(\bar z_t=(\gamma_t,\zeta_t)\in Z\), where \(\gamma_t=(\bar z_t)_0\in[0,1]\) is the coefficient on the level term \(f(x)-l\). The next lemma shows that the smooth CGO output satisfies the generic outer-loop requirements \eqref{eq:c1}--\eqref{eq:c3}.

\begin{lem}
 When CGO is applied to \eqref{cvxctr-model3} and LCG does not terminate at iteration $k$, the quantities \((\gamma_k, L_k, U_k)\) defined from its output \((x_k,\bar z_k,L_k,U_k)\) satisfy \eqref{eq:c1}-\eqref{eq:c3}.
\label{lem-oc}
\end{lem}
\proof{Proof.}
Fix a level estimate \(l\) and apply CGO to the specialized saddle problem \eqref{cvxctr-model3}, for which \(\bar f=0\), \(\bar h(x;l)=(f(x)-l,h(x))\), and \(\bar z_t=(\gamma_t,\zeta_t)\in Z\).
The scalar lower model built by CGO is $\underline h_0(x;l)= r_{0,0}\big[\ell_f(x_0,x)-l\big] + \sum_{i=1}^m r_{0,i}\ell_{h_i}(x_0,x),$ and, for \(t\ge 1\),
\begin{equation*}
\underline h_t(x;l)= (1-\alpha_t)\underline h_{t-1}(x;l) + \alpha_t\Big[r_{t,0}\big(\ell_f(x_{t-1},x)-l\big) + \sum_{i=1}^m r_{t,i}\ell_{h_i}(x_{t-1},x)\Big].
\end{equation*}	
Define $G_0(x):=r_{0,0}\ell_f(x_0,x)+\sum_{i=1}^m r_{0,i}\ell_{h_i}(x_0,x),$ and $\gamma_0:=r_{0,0},$ and, for \(t\ge 1\),
\[
G_t(x):=(1-\alpha_t)G_{t-1}(x)+\alpha_t\Big[r_{t,0}\ell_f(x_{t-1},x)+\sum_{i=1}^m r_{t,i}\ell_{h_i}(x_{t-1},x)\Big],
\]
\[
\gamma_t:=(1-\alpha_t)\gamma_{t-1}+\alpha_t r_{t,0}.
\]
An induction on \(t\) yields $\underline h_t(x;l)=G_t(x)-\gamma_t l,$ hence $\mathcal L_t(l):=\min_{x\in X}\underline h_t(x;l)= -\gamma_t l + \min_{x\in X}G_t(x),$ so \(\mathcal L_t(\cdot)\) is affine. Next, convexity of \(f\) and \(h_i\) gives \(\ell_f(x_{t-1},x)\le f(x)\) and \(\ell_{h_i}(x_{t-1},x)\le h_i(x)\). Therefore
\begin{equation*}
	\begin{aligned}
	r_{t,0}\big(\ell_f(x_{t-1},x)-l\big)+\sum_{i=1}^m r_{t,i}\ell_{h_i}(x_{t-1},x)
&\le
r_{t,0}\big(f(x)-l\big)+\sum_{i=1}^m r_{t,i}h_i(x)\\
&=
\left\langle (f(x)-l,h(x)),r_t\right\rangle.
	\end{aligned}
\end{equation*}
	Combining this with the recursion for \(\underline h_t\) and the update \(\bar z_t=(1-\alpha_t)\bar z_{t-1}+\alpha_t r_t\), another induction gives \(\underline h_t(x;l)\le \left\langle (f(x)-l,h(x)),\bar z_t\right\rangle \le \max_{(\gamma,z)\in Z}\ \gamma[f(x)-l]+\langle h(x),z\rangle\).
Minimizing over \(x\in X\) yields \(\mathcal L_t(l)\le \phi(l)\), which is \eqref{eq:c3}. Setting \(l=l_k\) gives \(L_k=\mathcal L_k(l_k)\le \phi(l_k)\), while \(\phi(l_k)\le U_k\) follows directly from the definition of \(U_k\) by evaluating the minimum in \(\phi(l_k)\) at \(x_k\). This proves \eqref{eq:c2}.

	Finally, if LCG has not terminated at outer iteration \(k\), then \(L_k>0\) by Lemma~\ref{lem5}. If \(\gamma_k=0\), then \(\underline h_k(x;l_k)=G_k(x)\) and \(\bar z_k=(0,\zeta_k)\in Z\). The previous bound gives \(G_k(x)=\underline h_k(x;l_k)\le \langle h(x),\zeta_k\rangle\) for all \(x\in X\). Let \(x^*\) be any feasible point of \eqref{cvxctr-model}. Since \(h_i(x^*)\le 0\) and \(\zeta_k\ge 0\), \(L_k=\min_{x\in X}G_k(x)\le G_k(x^*)\le \langle h(x^*),\zeta_k\rangle\le 0\), contradicting \(L_k>0\). Hence \(\gamma_k>0\), proving \eqref{eq:c1}. \Halmos
\endproof
We are now ready to establish the overall iteration complexity of 
	the LCG method.
	\begin{thm}\label{overall-thm-lcg}
	Under Assumptions \ref{asm:lcg-problem} and \ref{asm:cgo-smooth}, suppose that at each outer iteration the smooth CGO is run on \eqref{cvxctr-model3} with parameters \eqref{thm1-params1} until \(U_t-L_t\le (1-\mu)\epsilon\). Then the total number of CGO iterations required to find an \(\epsilon\)-solution \(\bar x\in X\) of \eqref{cvxctr-model} is bounded by \(\mathcal{O}\left(\tfrac{1}{\epsilon^2}\log(\tfrac{1}{\epsilon})\right)\).
	\end{thm}
\proof{Proof.}
By Lemma \ref{lem-oc}, the smooth CGO output satisfies \eqref{eq:c1}--\eqref{eq:c3}, so Theorem \ref{prop2} applies to the outer loop. Moreover, the specialization of Theorem \ref{thm:cgo-smooth-main} to \eqref{cvxctr-model3}, given in \eqref{eq:lcg-smooth-inner-bound}, shows that each inner CGO call needs \(\mathcal{O}(\epsilon^{-2})\) iterations to attain \(U_t-L_t\le (1-\mu)\epsilon\). Theorem \ref{prop2} then implies that \(\mathcal{O}(\log(1/\epsilon))\) outer iterations suffice to obtain \(U_k\le \epsilon\), and hence an \(\epsilon\)-solution of \eqref{cvxctr-model}. Combining the inner and outer bounds yields the overall complexity \(\mathcal{O}\left(\tfrac{1}{\epsilon^2}\log(\tfrac{1}{\epsilon})\right)\). \Halmos
\endproof

\subsection{Extension to Structured Nonsmooth Functions}
\label{sec-cgo-nonsmooth}

We now extend the same CGO template to structured nonsmooth saddle problems. The key idea is to smooth each nonsmooth component, run the same primal--dual/LMO framework on the smoothed model, and then bound the difference between the resulting smoothed certificates and the original nonsmooth ones.
The following assumption records the max-structure needed for this smoothing step.
\begin{assumption}\label{asm:cgo-nonsmooth}
For problem \eqref{model}, each functional has the max-structure
\begin{align}
\bar{f}(x) &= \max\limits_{y\in Y_0} \left\{\langle B_0x, y\rangle - \hat{f}(y)\right\}, \label{eq:ns-structure-f}\\
\bar{h}_i(x) &= \max\limits_{y\in Y_i} \left\{\langle B_ix, y\rangle - \hat{h}_i(y)\right\},\ i=1,\ldots,\bar m, \label{eq:ns-structure-h}
\end{align}
where each \(Y_i\) is closed convex, \(\hat f,\hat h_i\) are convex, and \(B_i\) are linear maps.
\end{assumption}
For each \(i=0,\ldots,\bar m\), let \(u_i:Y_i\to\mathbb{R}\) be a \(1\)-strongly convex distance generating function and define \(U_i(y):=u_i(y)-u_i(y_{u_i})-\langle \nabla u_i(y_{u_i}),y-y_{u_i}\rangle\), where \(y_{u_i}:=\argmin_{y\in Y_i}u_i(y)\).
We then apply Nesterov smoothing:
\begin{align}
\bar{f}_{\eta_0}(x) &:= \max\limits_{y\in Y_0}\ \left\{\langle B_0x, y\rangle - \hat{f}(y) - \eta_0U_0(y)\right\}, \label{eq:ns-smoothing-f}\\
\bar{h}_{i,\eta_i}(x) &:= \max\limits_{y\in Y_i}\ \left\{\langle B_ix, y\rangle - \hat{h}_i(y) - \eta_iU_i(y)\right\},\ i=1,\ldots,\bar m.\label{eq:ns-smoothing-h}
\end{align}
At iteration \(t\), Algorithm \ref{alg:CG oracle} is modified by replacing \(\ell_{\bar f}\) and \(\ell_{\bar h}\) with \(\ell_{\bar f_{\eta_0^t}}\) and \(\ell_{\bar h_{\eta^t}}\), and by using the extrapolated smoothed constraint model \(\tilde h_t=u_{t-1}+\lambda_t(u_{t-1}-u_{t-2})\) with \(u_{t-1}:=\ell_{\bar h_{\eta^{t-1}}}(x_{t-2},p_{t-1})\). This gives Algorithm \ref{alg:CG oracle nonsmooth} in Appendix \ref{sec:nonsmooth-convex}. The lower-bounding recursions keep the same form as \eqref{alg-cg-lowerf}--\eqref{alg-cg-lowerh}, now with smoothed linearizations.

The next proposition shows that this extension preserves the same order of inner-loop complexity.
\begin{prop}\label{prop:cgo-nonsmooth-main}
Under Assumption \ref{asm:cgo-nonsmooth}, let \(\alpha_t=2/(t+1)\), \(\lambda_t=(t-1)/t\), \(\tau_t=9\sqrt{t}\,\bar M_{B,U}D_{\bar X}\), and \(\eta_i^t=\|B_i\|D_{\bar X}/(\sqrt{t}D_{U_i})\), \(i=0,\ldots,\bar m\),
where \(D_{U_i}:=\left(\max_{y\in Y_i}U_i(y)\right)^{1/2}\) and
\(\bar{M}_{B,U}:=\left(\sum_{i=1}^{\bar m}M^2_{B_i,U_i}\right)^{1/2}\) with
\(M_{B_i,U_i}:=\|B_i\|(\|y_{u_i}\|+\sqrt{2}D_{U_i})\).
Then there exists a constant \(C_{\mathrm{ns}}>0\), independent of \(t\), such that \(U_t-L_t=\max_{w\in \bar X\times \bar Z}\bar Q_t(w_t,w)\le C_{\mathrm{ns}}/\sqrt{t+1}\) for \(t\ge 1\).
Consequently, to attain \(U_t-L_t\le (1-\mu)\epsilon\), CGO needs \(\mathcal{O}(\epsilon^{-2})\) iterations.
\end{prop}
The proof is given in Appendix~\ref{sec:nonsmooth-convex}.

\subsection{Overall Complexity in the Structured Nonsmooth Case}\label{sec:level-cg-overall-nonsmooth}
For the structured nonsmooth LCG subproblem, the same outer-loop argument still applies once the inner oracle satisfies \eqref{eq:c1}--\eqref{eq:c3}. Appendix \ref{sec:append-lcg} provides the nonsmooth analogue of Lemma \ref{lem-oc} together with a specialization of Theorem \ref{col-nonsmooth} to \eqref{cvxctr-model3}. Combining those results with Theorem \ref{prop2} yields the following counterpart of Theorem \ref{overall-thm-lcg}.
\begin{thm}\label{overall-thm-lcg-nonsmooth}
Under Assumptions \ref{asm:lcg-problem} and \ref{asm:cgo-nonsmooth}, suppose that at each outer iteration the structured nonsmooth variant of CGO is run on \eqref{cvxctr-model3} with the parameters of Proposition \ref{prop:cgo-nonsmooth-main} until \(U_t-L_t\le (1-\mu)\epsilon\). Then the total number of CGO iterations required to find an \(\epsilon\)-solution \(\bar x\in X\) of \eqref{cvxctr-model} is bounded by \(\mathcal{O}\left(\tfrac{1}{\epsilon^2}\log(\tfrac{1}{\epsilon})\right)\).
\end{thm}


\begin{rem}
It is worth mentioning here that as the output solution $x_k$ may not be a feasible solution such that $h(x_k) \le 0$, we develop a lower bound for $f(x_k) - f^*$. Under Slater's condition for \eqref{cvxctr-model}, i.e., there exists \(\tilde x\in \operatorname{ri}(X)\) such that \(h_i(\tilde x)<0\) for all \(i=1,\ldots,m\), such a bound is presented in Lemma \ref{lem-lowerbound} in Appendix \ref{sec:append-aux-lemma}. 
\end{rem}

\section{Nonconvex Conditional Gradient Methods}
\label{sec-nonconvex}
In this section, we consider the nonconvex functional constrained problem \eqref{cvxctr-model}, where \(f\) is nonconvex and \(h_i\), \(i=1,\cdots,m\), are convex. Since global or even local solution of this problem is generally intractable, our goal is instead to compute approximate first-order stationary points. To that end, we introduce an Inexact Proximal Point Level Conditional Gradient (IPP-LCG) method. The basic strategy is to reduce the nonconvex problem to a sequence of \emph{convex} functionally constrained subproblems, each of which can be handled by the projection-free LCG machinery from Section \ref{sec-level-cg}. In this way, the nonconvex method inherits the same LMO-based inner solves and remains projection-free throughout.

We first state the regularity conditions used throughout the nonconvex analysis.
\begin{assumption}\label{asm:nonconvex-smooth}
Function \(f\) is \(\underline L_f\)-smooth and \(M_f\)-Lipschitz continuous on \(X\). Each \(h_i\) is \(L_{h_i}\)-smooth and \(M_{h_i}\)-Lipschitz continuous on \(X\). The norm \(\|\cdot\|\) is induced by an inner product.
\end{assumption}
By \citep[Lemma 1.2.3]{nesterov2013introductory}, \(\underline L_f\)-smoothness of \(f\) implies the weakly convex property: for any \(x,y\in X\),
\begin{equation}
			f(x) -f(y) - \langle \nabla f(y), x- y\rangle \ge -\frac{\underline{L}_f}{2}\| x - y\|^2.
			\label{lower-curv-f}
\end{equation} 
Thus adding the quadratic term \(\underline L_f\|x-x_{j-1}\|^2\) convexifies the objective in each proximal subproblem.

\subsection{Inexact Proximal Point Method}

The IPP-LCG method applies LCG inexactly to a sequence of proximal convexifications of the original problem. Compared with classical proximal-point methods \citep{guler1992new,bertsekas2015convex,scutari2016parallel,lan2019accelerated,kong2019complexity,ma2019proximally,boob2020feasible}, our setting retains the general functional constraints and uses only projection-free inner solves. Given the current iterate \(x_{j-1}\in X\), we form the convex subproblem
\begin{equation}
	\begin{aligned}
		&\min\ f(x;x_{j-1}) := f(x) + \underline{L}_f \|x- x_{j-1}\|^2 \\
		&\st\ h_i(x) \le 0,\ i=1,\cdots,m,\\
		&\quad\quad  x\in X.
	\end{aligned}
	\label{convex-sub}
	\end{equation}
	To ensure uniform control of the dual variables for these subproblems, we additionally impose Slater's condition.
\begin{assumption}\label{asm:nonconvex-slater}
The Slater condition holds for \eqref{convex-sub}: there exists \(\tilde{x}\in \operatorname{ri}(X)\) such that \(h_i(\tilde{x})<0\) for all \(i=1,\dots,m\).
\end{assumption}
We use \((x^*_j, y^*_j)\) to denote a primal-dual optimal pair of \eqref{convex-sub}. At outer iteration \(j\), Algorithm \ref{alg:inexact proximal} calls LCG and returns an \((\tilde\epsilon,\tilde\epsilon)\)-optimal solution \(x_j\) of \eqref{convex-sub}, meaning that
\[
f(x_j;x_{j-1}) - f(x^*_j;x_{j-1})\le \tilde{\epsilon},
\qquad
\max_{i\in[1:m]}h_{i}(x_j)\le \tilde{\epsilon}.
\]
Among the candidates \(x_1,\ldots,x_J\), the method outputs an index \(\hat j \in \argmin_{j=1,\ldots,J}\{f(x_{j-1})-f(x_j)\}\). This choice is deliberate. Since \(\hat j\) minimizes the decrease \(f(x_{j-1})-f(x_j)\), we have \(f(x_{\hat j-1})-f(x_{\hat j}) \le J^{-1}\sum_{j=1}^J [f(x_{j-1})-f(x_j)] = [f(x_0)-f(x_J)]/J\).
Thus the selected iterate attains a decrease no larger than the average decrease across the \(J\) proximal subproblems, and this is the quantity used later to derive the stationarity bound in Theorem \ref{thm8}.

\begin{algorithm}[h]
	\caption{Inexact Proximal Point Level Conditional Gradient Method (IPP-LCG)} 
	\nonumber
		\begin{algorithmic}[1]
		\State {Initialization: $x_0\in X$.}
		\For {$j =1,2,\ldots, J$}
		\State {Call LCG to solve \eqref{convex-sub} and return a $(\tilde{\epsilon},\tilde{\epsilon})$-optimal solution $x_j$.}
		\EndFor
		\State{Select $\hat{j}$ such that $\hat j \in \argmin\limits_{j=1,\cdots,J} \{f(x_{ j -1}) - f(x_{ j})\}$.}
		\State {Terminate and return $x_{\hat{j}}$.}
	\end{algorithmic} 
	\label{alg:inexact proximal}
\end{algorithm}
We measure progress by the following near-KKT criterion.
\begin{defi}
For problem \eqref{cvxctr-model},
\begin{enumerate}
    \item[(i)] 
	    $x'$ is an $(\epsilon, \delta)$-near-KKT point if $x'\in X$ and there exists $(x, y)$ such that $h_i(x) \le 0,x\in X, y_i \ge 0, i=1,\cdots, m$ and 
	\begin{equation}
	\sum\limits_{i=1}^{m}|y_ih_i(x)|\le \epsilon,\quad \left[d\left(\nabla f(x)+\sum\limits_{i=1}^m y_i\nabla h_i(x),-N_X(x)\right)\right]^2\le \epsilon,\quad \|x'-x\|^2\le \delta.
	\label{kkt-criteria}
	\end{equation}
	where $\epsilon,\delta >0$, $d(\cdot,\cdot)$ denotes the distance between two sets $A$ and $B$ such that $d(A,B) := \min\limits_{a\in A, b\in B} \|a - b\|$.
	\item[(ii)] $x$ is an $\epsilon$-KKT point (paired with $y$) if it satisfies the first two criteria in \eqref{kkt-criteria} with $h_i(x)\le 0, x\in X, y_i\ge 0, i =1,\cdots,m$.
	\end{enumerate}
		\label{def-inexact-kkt-point}
\end{defi}
To pass from inexact solutions of the convex subproblems to a near-KKT guarantee for the original nonconvex problem, we need uniform control of the multipliers associated with those subproblems. The next lemma provides exactly this bound under Slater's condition.
\begin{lem}\label{lem7}\label{lem: bounded-dual}
Under Assumption \ref{asm:nonconvex-slater}, define the Slater margin \(\sigma := -\max_{i\in[m]} h_i(\tilde{x}) > 0\), let \(\overline{V}:= \max_{x,y\in X} \|x-y\|^2\), and set \(B :=  [f(\tilde{x}) - f^* + \underline{L}_{f}\overline{V}]/\sigma\). Then every dual optimal solution \(y_j^*\) of \eqref{convex-sub} satisfies \(\|y_j^*\|_1 \le B\).
\end{lem}
The proof is given in Appendix \ref{append-sub-sec-IPP-LCG}.
We are now ready to state the main convergence guarantee for IPP-LCG.
\begin{thm}
Under Assumptions \ref{asm:nonconvex-smooth} and \ref{asm:nonconvex-slater}, let \(D_X:=\max_{x,y\in X}\|x-y\|\).
If \(J\ge 16\underline L_fM_fD_X/\bar\epsilon\) and \(\tilde\epsilon\le \min\{\bar\epsilon/[2(1+B)],\underline L_f\bar\delta/[2(1+B)]\}\) in Algorithm \ref{alg:inexact proximal}, then the total number of CGO iterations performed by 
	the IPP-LCG method to compute an $(\bar{\epsilon},\bar{\delta})$-near-KKT point of problem \eqref{cvxctr-model} is bounded by $\mathcal{O}\left(\bar{\epsilon}^{-1}\cdot\max\left\{\bar{\delta}^{-2}\log(\bar{\delta}^{-1}),\bar{\epsilon}^{-2}\log(\bar{\epsilon}^{-1})\right\}\right)$.
\label{thm8}
\end{thm}
The proof is given in Appendix \ref{append-sub-sec-IPP-LCG}.
Theorem \ref{thm8} should be read as follows: the algorithm returns \(x_{\hat j}\), and the analysis guarantees that this point lies within distance at most \(\sqrt{\bar\delta}\) of a point \(x_{\hat j}^*\) satisfying the corresponding \(\bar\epsilon\)-KKT conditions.
The following corollary records the common choice in which both tolerances are set to the same target accuracy.
\begin{cor}
Under the assumptions of Theorem \ref{thm8}, if we set \(\bar{\epsilon}=\bar{\delta}=\epsilon\), then IPP-LCG computes an \((\epsilon,\epsilon)\)-near-KKT point using a total number of CGO iterations bounded by \(\mathcal{O}\left(\epsilon^{-3}\log(\epsilon^{-1})\right)\).
\label{cor:ipp-lcg-eps}
\end{cor}
The proof is given in Appendix \ref{append-sub-sec-IPP-LCG}.
\begin{rem}\label{rmk1}
For \(X=\mathbb{R}^n\), Theorem \ref{thm8} implies that the returned point \(x_{\hat j}\) itself has small first-order residuals: its stationarity residual, complementarity residual, and feasibility violation are all small. In this sense, although Theorem \ref{thm8} is stated through a nearby \(\bar\epsilon\)-KKT point \(x_{\hat j}^*\), the actual output \(x_{\hat j}\) also behaves as an approximate stationary point in the unconstrained-set case. Let \((x_{\hat j}^*,y_{\hat j}^*)\) be the point guaranteed by Theorem \ref{thm8}. Then \(\|x_{\hat j}-x_{\hat j}^*\|^2\le \bar\delta\), \(x_{\hat j}^*\) is a \(\bar\epsilon\)-KKT point, and \(\|y_{\hat j}^*\|_1\le B\). Since \(N_X(x)=\{0\}\) for \(X=\mathbb{R}^n\), Lipschitz continuity of \(\nabla f\) and \(\nabla h_i\) implies
\(\|\nabla f(x_{\hat j})+\sum_{i=1}^m y_{\hat j,i}^*\nabla h_i(x_{\hat j})\|^2=O(\bar\delta+\bar\epsilon)\), \(\sum_{i=1}^m |y_{\hat j,i}^* h_i(x_{\hat j})|=O(\sqrt{\bar\delta}+\bar\epsilon)\), and \(\max_{i\in[m]}[h_i(x_{\hat j})]_+=O(\sqrt{\bar\delta})\).
The proof is given in Appendix \ref{append-sub-sec-IPP-LCG}.
\end{rem}
It is worth noting that IPP-LCG can also be generalized to solve structured nonsmooth problems by applying LCG to the resulting nonsmooth convex subproblems. The convergence analysis in this case is essentially the same as that shown in Theorem \ref{thm8}.

\begin{rem}
In practice, it is often difficult to know the lower curvature parameter \(\underline L_f\) exactly. One may therefore treat \(\underline L_f\) as a tuning parameter or use conservative estimates, for example by increasing it until the proximal subproblems behave stably. When second-order information is available, one may also use curvature estimates based on Hessian information, such as the smallest eigenvalue of \(\nabla^2 f\), as a heuristic proxy for \(\underline L_f\); see, e.g., \citet{liang2021average}.
\end{rem}

\section{Applications and Numerical Experiments}\label{sec:applications}\label{sec:numerical}
This section instantiates the proposed framework in portfolio selection and intensity-modulated radiation therapy (IMRT) and evaluates the practical behavior of LCG and IPP-LCG in both settings. In portfolio selection, the experiments examine whether the proposed methods can attain favorable out-of-sample risk under explicit sparsity targets on convex and smooth nonconvex training models. In IMRT, the experiments examine whether the methods can construct sparse treatment plans with good clinical feasibility. An anonymous repository containing the implementation and instructions for reproducing the experiments is available \href{https://anonymous.4open.science/r/projection-free-functional-constraints-C683/README.md}{here}. All experiments were run in Python 3.9.18 on Red Hat Enterprise Linux 9.4 on a server with two AMD EPYC 9334 CPUs and 377\,GiB RAM.

\subsection{Portfolio Selection}\label{sec:app-portf}\label{sec:numerical-portf}
We study a risk-averse portfolio problem in which the investor seeks to control benchmark-underperformance risk while limiting the number of selected assets. The numerical study considers both a convex, \CVaR-based surrogate and a smooth nonconvex surrogate of the underlying step-risk objective, together with a cap-type convex constraint that encodes a target support level \(\Psi\).

\subsubsection{Models}\label{subsubsec:portf-model}
Consider a long-only portfolio of \(N\) risky assets with random returns \(r_i\) and random benchmark return \(R\). Let \(x_i\) denote the weight of asset \(i\), let \(x=(x_1,\ldots,x_N)\), and let
\(
X=\{x\ge 0:\ \sum_{i=1}^N x_i = 1\}
\)
be the simplex of feasible portfolio weights. We measure underperformance through the step-risk
\begin{equation}
\mathrm{Risk}(x)=\mathbb{P}\!\left(R-\sum_{i=1}^N r_i x_i>0\right)
=\mathbb{E}\!\left[\mathbb{I}\!\left\{R-\sum_{i=1}^N r_i x_i>0\right\}\right],\label{eq_risk_def}
\end{equation}
namely, the probability that the portfolio return falls below the benchmark. Given \(K\) historical samples \(\{(R_k,r_{1k},\ldots,r_{Nk})\}_{k=1}^K\), we consider the sample-average approximation with a target support size \(\Psi\). Specifically, letting \(s(x):=|\{i:\ x_i>0\}|=\|x\|_0\), we start from
\begin{equation}
\min_{x\in X,\ s(x)\le \Psi}\ f(x):=\frac{1}{K}\sum_{k=1}^K \mathbb{I}{\Big\{R_k - \sum_{i=1}^N r_{ik} x_i > 0\Big\}}.
\label{portf-var-sa}
\end{equation}

The exact cardinality budget \(s(x)\le \Psi\) is combinatorial and therefore lies outside the convex functional-constraint template of Section~\ref{sec-level-cg}. In the experiments, we replace it by a tractable convex proxy over the augmented variable \((x,v)\in X\times[\underline v,\bar v]\):
\begin{equation}
g_{\Psi}(x,v):=Nv+\frac{1}{\Psi}\sum_{i=1}^N [x_i-v]_+-\frac{N}{\Psi}\le 0,
\qquad
\underline v\le v\le \bar v,
\label{portf-sparsity-surrogate}
\end{equation}
where \([t]_+:=\max\{t,0\}\). Here \(v\) plays the role of a common threshold, and the hinge terms \([x_i-v]_+\) measure how much each weight exceeds that threshold. This does not reproduce the exact support constraint, but it discourages weight from being spread across many assets while remaining convex in \((x,v)\), since each hinge term is convex and the remaining terms are affine. Throughout the portfolio experiments, we set \(\underline v=10^{-4}\) and \(\bar v=1/\Psi\). The offset \(N/\Psi\) calibrates \eqref{portf-sparsity-surrogate} so that the equal-weight \(\Psi\)-asset portfolio lies on the surrogate boundary. In the reported results, however, sparsity is always evaluated by the realized support \(s(x)\), not by the surrogate value.

\paragraph{Convex and smooth relaxations.}
The step-risk objective in \eqref{portf-var-sa} is discontinuous and nonconvex, so we study two surrogate training models (details in Appendix~\ref{appen-num-portf}).

\emph{Smooth nonconvex surrogate.} For a smoothing parameter \(\theta>0\), define \(\tilde{f}_{\theta}(x):= \tfrac{1}{K}{\textstyle \sum}_{k=1}^K \sigma\!\left(\tfrac{R_k-{\textstyle \sum}_{i=1}^N r_{ik}x_i}{\theta}\right)\) with \(\sigma(t):=\tfrac{1}{1+\exp(-t)}\). As \(\theta\to 0\), \(\tilde f_\theta\) approximates the step function \(f\) pointwise. We use the following smooth but nonconvex risk surrogate in the portfolio training models:
\begin{equation}
\min_{x\in X,\ v\in[\underline v,\bar v]}\ \tilde f_\theta(x)\quad \text{s.t.}\quad g_{\Psi}(x,v)\le 0.
\label{portf-nonconvex}
\end{equation}

\emph{Convex nonsmooth surrogate.} For a confidence level \(\alpha\in(0,1)\), we replace the indicator-based objective by the empirical \CVaR\ of the underperformance loss \(R_k-r_k^\top x\), which yields the convex and nonsmooth model
\begin{equation}
\min_{x\in X,\ u\in[\underline u,\bar u],\ v\in[\underline v,\bar v]}\ \left\{u+\frac{1}{\alpha K}\sum_{k=1}^K \bigg[R_k-\sum_{i=1}^N r_{ik}x_i-u\bigg]_+\right\}\quad \text{s.t.}\quad g_{\Psi}(x,v)\le 0,
\label{portf-convex}
\end{equation}
where \([\underline u,\bar u]\) is a data-dependent interval used to keep the feasible set compact. A brief derivation of \eqref{portf-convex} as a convex surrogate of \eqref{portf-var-sa}, together with the convexity argument, is given in Appendix~\ref{appen-num-portf}. Hence \eqref{portf-convex} fits the convex functional-constrained template \eqref{cvxctr-model}, and its objective falls under the structured nonsmooth setting in Section~\ref{sec-cgo-nonsmooth}.

\paragraph{Training formulations used in the experiments.}
In the numerical study, we use two training problems: the convex \CVaR\ formulation \eqref{portf-convex}\phantomsection\label{model1} and the smooth nonconvex formulation \eqref{portf-nonconvex}\phantomsection\label{model2}.

\subsubsection{Tests on Stock Market Dataset}
We test the above portfolio models on historical weekly returns from six major stock markets (from \textit{Thomson Reuters Datastream} and the \textit{Fama \& French Data Library}). Each dataset provides weekly asset returns $\{r_{ik}\}$ and an index return series $\{R_k\}$ over $K$ weeks. Since the target cardinality $\Psi$ is not specified in the raw data, we set it using a standard rule-of-thumb that scales with the universe size:
\(\Psi = \lfloor 0.2\,N \rfloor\) if \(N \le 100\) and \(\Psi = \lfloor 0.05\,N \rfloor\) otherwise.
Appendix Table \ref{tab:stock-data} summarizes dataset characteristics (see also \citep{bruni2016real}).

We report the in-sample step-risk \eqref{portf-var-sa} and out-of-sample performance on a single 70/30 train/test split. For each dataset, all methods are evaluated on the same split. Throughout, the out-of-sample metric is
\(\mathrm{StepRisk}_{\delta}(x):=\tfrac{1}{K_{\mathrm{test}}}{\textstyle \sum}_{k\in \mathrm{test}}\mathbb{I}\{R_k-{\textstyle \sum}_{i=1}^N r_{ik}x_i>\delta\}\),
which estimates the probability of underperforming the index by more than \(\delta\). Unless stated otherwise, \(\delta=0.0005\). Appendix Table~\ref{tab:portf-ref-risk} reports reference values for simple benchmark portfolios. For sparsity reporting we use \(s(x):=\big|\{i:\ x_i>10^{-4}\}\big|\), i.e., weights below \(10^{-4}\) are treated as numerical dust.

\paragraph{Evaluation protocol (time-budgeted sparse solver comparison).}
For the convex \CVaR\ model, we compare LCG with CoexDurCG \citep{lan2020conditional}, FISTA, PGD, PGD-IHT, mirror descent with entropy prox (MD(ent)), Greedy+Refit, IRL1, and PenPGD. For the smooth nonconvex model, we compare IPP-LCG with PGD, PGD-IHT, MD(ent), Greedy+Refit, IRL1, and PenPGD; Appendix~\ref{sec:append:portf-baselines} describes these baseline methods in more detail. We compare all methods under a common wall-clock budget of 5 seconds. For the projection-free methods, the training problems optimize over \((x,v)\) and use the convex surrogate constraint \(g_{\Psi}(x,v)\le 0\), while the target cardinality \(\Psi\) is used for evaluation across all methods. We select the best recorded iterate satisfying \(s(x)\le \Psi\) using training \(\mathrm{StepRisk}_{\delta}\) as the primary criterion; if a method never reaches the target, we penalize support violation. Appendix~\ref{append:portf-eval} gives the exact iterate-selection rule and Appendix Table~\ref{tab:portf-time-iters} reports the per-dataset update counts completed within the same 5-second budget.

\paragraph{Results under a fixed time budget.}
Tables \ref{tab:portf-time-convex} and \ref{tab:portf-time-nonconvex} report the out-of-sample step-risk at \(\delta=0.0005\), together with the support violation \( [s(x)-\Psi]_+ \) and achieved support size \(s(x)\) (\(\#\)ass.), after 5 seconds per method. Bold indicates the lowest out-of-sample \(\mathrm{StepRisk}_{\delta}\) among methods that satisfy the target(s).

\begin{table}[h]
\centering
\caption{Selecting the best iterate with $\|x\|_0\le \Psi$ (penalized if never reached), time-budget comparison on the convex portfolio model (5s per method). We use penalty $\rho=1$ when no iterate reaches $\Psi$. We report step-risk at threshold $\delta=0.0005$ and the support violation $[s(x)-\Psi]_+$. Bold indicates the best (lowest) StepRisk among methods with $[s(x)-\Psi]_+=0$ (infeasible methods are excluded from bolding; if no method is feasible, no entry is bolded).}
\resizebox{\columnwidth}{!}{%
\begin{tabular}{l|ccccccccccccccc}
\toprule
 & \multicolumn{3}{c}{LCG} & \multicolumn{3}{c}{CoexDurCG} & \multicolumn{3}{c}{FISTA} & \multicolumn{3}{c}{Greedy+Refit} & \multicolumn{3}{c}{IRL1} \\
\midrule
Instance & StepRisk$_{\delta}$(test) & $[s(x)-\Psi]_+$ & \# ass. & StepRisk$_{\delta}$(test) & $[s(x)-\Psi]_+$ & \# ass. & StepRisk$_{\delta}$(test) & $[s(x)-\Psi]_+$ & \# ass. & StepRisk$_{\delta}$(test) & $[s(x)-\Psi]_+$ & \# ass. & StepRisk$_{\delta}$(test) & $[s(x)-\Psi]_+$ & \# ass. \\
\midrule
DJ & 0.4499 & 0 & 5 & \textbf{0.4377} & 0 & 5 & 0.4523 & 23 & 28 & 0.4792 & 0 & 5 & 0.4817 & 1 & 6 \\
FF49 & \textbf{0.4413} & 0 & 9 & 0.4642 & 0 & 6 & 0.4140 & 20 & 29 & 0.4670 & 0 & 9 & 0.4685 & 0 & 8 \\
ND100 & \textbf{0.3855} & 0 & 14 & 0.4078 & 0 & 11 & 0.3687 & 34 & 50 & 0.4860 & 0 & 16 & 0.4134 & 0 & 15 \\
FTSE100 & \textbf{0.3935} & 0 & 16 & 0.4120 & 0 & 15 & 0.3889 & 29 & 45 & 0.4352 & 0 & 11 & 0.5324 & 0 & 12 \\
SP500 & \textbf{0.3911} & 0 & 18 & 0.4190 & 0 & 14 & 0.3631 & 65 & 87 & 0.5642 & 0 & 12 & 0.4581 & 0 & 15 \\
NDComp & \textbf{0.4466} & 0 & 19 & 0.4660 & 0 & 19 & 0.4369 & 993 & 1053 & 0.4951 & 0 & 49 & 0.5097 & 0 & 58 \\
\bottomrule
\end{tabular}%
}
\vspace{0.5ex}
\resizebox{0.9\columnwidth}{!}{%
\begin{tabular}{l|cccccccccccc}
\toprule
 & \multicolumn{3}{c}{PenPGD} & \multicolumn{3}{c}{PGD} & \multicolumn{3}{c}{PGD-IHT} & \multicolumn{3}{c}{MD(ent)} \\
\midrule
Instance & StepRisk$_{\delta}$(test) & $[s(x)-\Psi]_+$ & \# ass. & StepRisk$_{\delta}$(test) & $[s(x)-\Psi]_+$ & \# ass. & StepRisk$_{\delta}$(test) & $[s(x)-\Psi]_+$ & \# ass. & StepRisk$_{\delta}$(test) & $[s(x)-\Psi]_+$ & \# ass. \\
\midrule
DJ & 0.5086 & 2 & 7 & 0.4181 & 21 & 26 & 0.5012 & 0 & 5 & 0.4059 & 23 & 28 \\
FF49 & 0.4685 & 0 & 8 & 0.0029 & 40 & 49 & 0.5086 & 0 & 9 & 0.0029 & 40 & 49 \\
ND100 & 0.4134 & 0 & 15 & 0.3743 & 42 & 58 & 0.4749 & 0 & 16 & 0.3631 & 66 & 82 \\
FTSE100 & 0.5324 & 0 & 12 & 0.4398 & 40 & 56 & 0.5139 & 0 & 16 & 0.4213 & 64 & 80 \\
SP500 & 0.4302 & 0 & 16 & 0.3631 & 168 & 190 & 0.4749 & 0 & 22 & 0.3017 & 409 & 431 \\
NDComp & 0.5097 & 0 & 58 & 0.5340 & 161 & 221 & 0.5243 & 0 & 60 & 0.4466 & 1143 & 1203 \\
\bottomrule
\end{tabular}%
}
\label{tab:portf-time-convex}
\end{table}

\begin{table}[ht]
\centering
\caption{Selecting the best iterate with $\|x\|_0\le \Psi$ (penalized if never reached), time-budget comparison on the nonconvex portfolio model (5s per method). We use penalty $\rho=1$ when no iterate reaches $\Psi$. We report step-risk at threshold $\delta=0.0005$ and the support violation $[s(x)-\Psi]_+$. Bold indicates the best (lowest) StepRisk among methods with $[s(x)-\Psi]_+=0$ (infeasible methods are excluded from bolding; if no method is feasible, no entry is bolded).}
\resizebox{\columnwidth}{!}{%
\begin{tabular}{l|cccccccccccc}
\toprule
 & \multicolumn{3}{c}{IPP-LCG} & \multicolumn{3}{c}{Greedy+Refit} & \multicolumn{3}{c}{IRL1} & \multicolumn{3}{c}{PenPGD} \\
\midrule
Instance & StepRisk$_{\delta}$(test) & $[s(x)-\Psi]_+$ & \# ass. & StepRisk$_{\delta}$(test) & $[s(x)-\Psi]_+$ & \# ass. & StepRisk$_{\delta}$(test) & $[s(x)-\Psi]_+$ & \# ass. & StepRisk$_{\delta}$(test) & $[s(x)-\Psi]_+$ & \# ass. \\
\midrule
DJ & \textbf{0.4377} & 0 & 5 & \textbf{0.4377} & 0 & 5 & 0.4034 & 14 & 19 & 0.3985 & 14 & 19 \\
FF49 & 0.4971 & 0 & 9 & 0.4656 & 0 & 9 & 0.4284 & 10 & 19 & 0.4284 & 14 & 23 \\
ND100 & \textbf{0.4022} & 0 & 16 & 0.4525 & 0 & 16 & 0.4022 & 18 & 34 & 0.3855 & 20 & 36 \\
FTSE100 & \textbf{0.3657} & 0 & 16 & 0.3981 & 0 & 15 & 0.3472 & 11 & 27 & 0.3472 & 12 & 28 \\
SP500 & 0.4413 & 0 & 2 & 0.4469 & 0 & 13 & 0.3743 & 28 & 50 & 0.3575 & 34 & 56 \\
NDComp & \textbf{0.4417} & 0 & 21 & 0.4466 & 0 & 23 & 0.5000 & 13 & 73 & 0.4903 & 33 & 93 \\
\bottomrule
\end{tabular}%
}
\vspace{0.5ex}
\resizebox{0.8\columnwidth}{!}{%
\begin{tabular}{l|ccccccccc}
\toprule
 & \multicolumn{3}{c}{PGD} & \multicolumn{3}{c}{PGD-IHT} & \multicolumn{3}{c}{MD(ent)} \\
\midrule
Instance & StepRisk$_{\delta}$(test) & $[s(x)-\Psi]_+$ & \# ass. & StepRisk$_{\delta}$(test) & $[s(x)-\Psi]_+$ & \# ass. & StepRisk$_{\delta}$(test) & $[s(x)-\Psi]_+$ & \# ass. \\
\midrule
DJ & 0.4205 & 1 & 6 & 0.4645 & 0 & 5 & 0.4034 & 23 & 28 \\
FF49 & 0.4470 & 2 & 11 & \textbf{0.4585} & 0 & 9 & 0.1203 & 40 & 49 \\
ND100 & 0.4637 & 0 & 16 & 0.4302 & 0 & 16 & 0.4078 & 66 & 82 \\
FTSE100 & 0.4398 & 0 & 15 & 0.4491 & 0 & 16 & 0.3889 & 67 & 83 \\
SP500 & \textbf{0.4246} & 0 & 19 & 0.4469 & 0 & 22 & 0.4022 & 420 & 442 \\
NDComp & 0.4660 & 0 & 44 & 0.5097 & 0 & 53 & 0.4515 & 1143 & 1203 \\
\bottomrule
\end{tabular}%
}
\label{tab:portf-time-nonconvex}
\end{table}

\paragraph{Takeaways.}
In the convex setting (Table \ref{tab:portf-time-convex}), LCG attains the best feasible out-of-sample \(\mathrm{StepRisk}_{\delta}\) on five of the six datasets, with CoexDurCG best on DJ, while consistently meeting the sparsity budget. In the nonconvex setting (Table \ref{tab:portf-time-nonconvex}), IPP-LCG attains the best feasible \(\mathrm{StepRisk}_{\delta}\) on four of the six datasets, while PGD/PGD-IHT are best on the remaining two.

\subsection{IMRT Treatment Planning}\label{sec:app-imrt}\label{sec:numerical-imrt}
In radiation therapy, a patient receives dose from a linear accelerator operating at angles \(a \in A\) and using apertures \(e \in E_a\). Treatment planning selects angles, apertures, and intensities to target tumor tissue while limiting overdose to healthy structures.

\subsubsection{Models}
A desirable treatment plan uses only a small number of angles. Following \citet{lan2020conditional}, we impose the group-sparsity constraint
\begin{equation}
\sum\limits_{a\in A} \max_{e\in E_a} y_{a,e} \le \Phi,
\label{group-sparse}
\end{equation}
where $y_{a,e}$ is the intensity of aperture $e$. Clinical criteria are imposed on voxelized structures; for example, ``PTV68: V68 $\ge 95\%$'' requires at least \(95\%\) of voxels in PTV68 to receive 68 Gy (gray, the standard unit of absorbed radiation dose) or more, while ``PTV68: V74.8 $\le 10\%$'' limits the fraction receiving at least 74.8 Gy to \(10\%\). Such criteria can be written as quantile-type underdose/overdose requirements.

\paragraph{Convex Formulation.}
We follow \citet{lan2020conditional} and use \CVaR\ approximations of the quantile criteria to obtain a convex model; see Appendix \ref{sec:appen-imrt-cvx}.

\paragraph{Nonconvex Formulation.}
To reduce approximation bias, we also consider a nonconvex model that uses smooth surrogates of the original quantile criteria in the objective together with the group-sparsity constraint \eqref{group-sparse}; see Appendix \ref{sec:appen-imrt-cvx}.

\subsubsection{Tests on Synthetic Dataset} \label{sec:numerical-imrt-synthetic}
We compare CoexDurCG and LCG on the convex formulation with \(\Phi = 0.005\) in \eqref{group-sparse} using the synthetic instances from Appendix \ref{sec:imrt-synthetic}. This isolates the effect of the proposed level-set framework within projection-free constrained methods. Table \ref{tab:ran-result-coex} reports the objective value, total constraint violation \(\|h(\hat{x})\|_2\), group-sparsity violation \(\|h_s\|_2\), clinical-constraint violation \(\|h_c\|_2\), and CPU time after 1000 iterations. The runtimes are similar, while LCG consistently returns smaller constraint violations, especially in the clinical constraints.
\begin{table}[h]
\caption{Results of applying CoexDurCG and LCG on synthetic data with $\Phi =0.005$ at iteration $1000$.	}
	\centering
	{\small
	\begin{tabular}{c|ccccc|ccccc|}
		\toprule
		 \multirow{2}{4em}{\centering{Instance}} &\multicolumn{5}{c|}{\centering{CoexDurCG}} &\multicolumn{5}{c|}{\centering{LCG}} \\
		{ }  &$f(\hat{x})$ &  $\|h(\hat{x})\|_2$ & $\|h_s\|_2$ &$\|h_c\|_2$ & Time (s) &$f(\hat{x})$ &  $\|h(\hat{x})\|_2$ & $\|h_s\|_2$ &$\|h_c\|_2$& Time (s)\\
		\midrule
		1   &0.0193  & 0.984   & 0.641 &0.747 & 926  &0.0193  & 0.528 &0.421 & 0.319 & 924 \\
		\midrule
		2   &0.0166  & 1.643 &0.614 &1.524  & 996 &0.019  & 0.763  & 0.402 &0.649& 908 \\
		\midrule
		3  & 0.0467  & 1.043  & 0.205 & 1.023 &4889  &0.047  & 0.476 &0.169 &0.445& 4834\\
	\midrule
		4  & 0.0465  & 3.193  & 0.208 & 3.186 &4867  &0.0435  & 0.984  & 0.175 &0.968 & 4871\\
		\bottomrule
	\end{tabular}}
	\label{tab:ran-result-coex}
    \vspace{-0.5cm}
\end{table}

\subsubsection{Tests on Prostate Dataset}
We next consider a large prostate-cancer instance in the treatment-planning setting of \citet{lan2020conditional}. It contains \(3{,}047{,}040\) voxels, \(180\) angles, and on average \(155\) beamlets per angle, so the implicit data matrices exceed \(3{,}047{,}040 \times 155 \times 180\). The instance includes 10 clinical criteria over six structures, with a narrow PTV68 dose window \([68,74.8]\) Gy and a strong preference for using few angles and at most 100 apertures. If an angle has an \(m\times n\) beamlet grid, then the number of potential apertures is \((\tfrac{n(n-1)}{2})^m\); with average \(m=n=10\), this yields an implicit variable space of size \(180\times 45^{10}\). As in column-generation IMRT optimization \citep{romeijn2005column,romeijn2008intensity}, we therefore do not enumerate apertures and instead solve a structured pricing subproblem inside each LMO step; see Appendix \ref{sec:appen-imrt-cvx}.

\paragraph{Results of the Convex Formulation.}
We apply LCG to the convex formulation. Table \ref{tab:prostate-sum} summarizes the resulting treatment plans through the number of iterations, selected angles, and apertures. As \(\Phi\) decreases, the plan uses fewer angles, which is the intended effect of the group-sparsity constraint.
\begin{table}[h]
\vspace{-0.3cm}
\caption{Treatment plans constructed by LCG on Prostate dataset with different $\Phi$.	}
	\centering
	{\small
	\begin{tabular}{cccccc}
		\toprule
		 {\centering{$\Phi$}} &{\centering{$\#$ of iter.}} &{\centering{$\#$ of angles}} &$\#$ of apertures\\
		\midrule
		1.0 & 100 & 27 & 99 \\
		0.5 & 85 & 17 & 84\\
		   0.005  &63  &6  & 62  \\
		   0.0005  &78  & 5  & 77  \\
		\bottomrule
\end{tabular}}
	\label{tab:prostate-sum}
    \vspace{-0.5cm}
\end{table}
Table \ref{tab:prostate-result} reports the dose-volume quantities for each clinical criterion under different sparsity parameters \(\Phi\). All entries in Table \ref{tab:prostate-result} (and Table~\ref{tab:prostate-noncvx-result-budget}) are percentages of voxels satisfying the stated dose event. Among the tested values, \(\Phi=0.005\) is the smallest sparsity parameter that still leaves only one violated listed criterion, namely the PTV68 overdose condition V\(74.8\le 10\%\). Compared with \(\Phi=1.0\) and \(\Phi=0.5\), it uses substantially fewer angles and apertures while also yielding the smallest remaining overdose violation among those three settings (Tables \ref{tab:prostate-sum} and \ref{tab:prostate-result}). Tightening further to \(\Phi=0.0005\) saves only one additional angle but introduces an extra target-coverage violation at PTV68 V68. From a clinical perspective, a limited violation of the overdose constraint on the target is often deemed acceptable in the presence of other, more important, properties of a treatment plan. This is due to the fact that this constraint aims to limit non-specific side effects of treating the target rather than the core goals of delivering sufficient dose to the target and limiting dose to specific critical structures in order to preserve their functionality. We therefore treat \(\Phi=0.005\) as the main near-feasible trade-off point in the convex study. 
\begin{table}[h]
\vspace{-0.3cm}
\caption{Results of applying LCG on Prostate dataset.	}
	\centering
	\resizebox{0.8\columnwidth}{!}{
	\begin{tabular}{ccccccccc}
		\toprule
 \multirow{4}{2em}{\centering{$\Phi$}}  &{\centering{PTV56}}  &{\centering{PTV68}} & Rectum &{\centering{Bladder}}  &{\centering{Lft. femoral head}} &{\centering{Rht. femoral head}}\\
		{ }  &V56$\ge 95\%$ &V68$\ge 95\%$ &V30$\le 80\%$ &V40$\le 70\%$   &V50$\le 1\%$  &V50$\le 1\%$ \\
		{ }  &{ } &V74.8$\le 10\%$ &V50$\le 50\%$ &V65$\le 30 \%$ &{} &{}\\
		{ }  &{ } &{} &V65$\le 25\%$ &{} &{} &{}\\
		\midrule
		\multirow{3}{2em}{\centering{1.0}}
		  & 99.97  & 96.47 & 68.25  &53.65 &0.11 & 0.10 \\
		  &  &15.93  &21.88  &21.79 &  &   \\
          &  &  &6.01 & & & \\
		\midrule
		\multirow{3}{2em}{\centering{0.5}}
		  & 99.96  & 95.36 & 71.26 &52.39  &1.00 & 0.23 \\
		  &  &14.23  &24.55  &21.15 &  &   \\
          &  &  &5.33 & & & \\
		\midrule
		\multirow{3}{2em}{\centering{0.005}}
		  & 99.87  & 95.44  & 77.78  &55.20  &0.24 & 0.00 \\
		  &  &12.63 & 37.86  &22.87 &  &   \\
          &  &  &5.61 & & & \\
		\midrule
\multirow{3}{2em}{\centering{0.0005}}
		  & 99.94  &90.56  & 79.98 &61.98 &0.03 & 0.22 \\
		  &   &11.51 & 38.04  &25.56 &  & \\
		  &   &  &10.71  &   &  \\
		\midrule
	\end{tabular}}
	\label{tab:prostate-result}
    \vspace{-0.5cm}
\end{table}

\paragraph{Results of the Nonconvex Formulation.}
We next study the nonconvex refinement phase on the prostate instance. Based on the convex comparison above, we use \(\Phi=0.005\) as the primary sparsity setting, since it is the sparsest tested convex plan with only one remaining listed violation. We then compare cold-start IPP-LCG and LCG-warm-start IPP-LCG under a common nonconvex wall-clock budget of 180 seconds. In the warm-started variant, the initialization is the convex LCG iterate obtained at iteration 63 for \(\Phi=0.005\) (Tables \ref{tab:prostate-sum} and \ref{tab:prostate-result}). For each protocol, we report the best iterate reached within the budget rather than the last iterate, so Tables \ref{tab:prostate-noncvx-sum-budget} and \ref{tab:prostate-noncvx-result-budget} compare best budgeted iterates. The budget applies only to the nonconvex phase, so the warm-started row does not include the time spent constructing the convex initialization. Because the budget is checked only after completed iterations, the reported nonconvex times can slightly exceed 180 seconds.
\begin{table}[h]
\caption{Best iterates under a fixed nonconvex wall-clock budget of 180.0s for IPP-LCG on the prostate dataset.}
\centering
{\small
\begin{tabular}{lccccccc}
\toprule
Type & $\Phi$ & $\#$ of iter. & $\#$ of angles & $\#$ of apertures & Time (s) & $\#$ viol. & Max viol.\\
\midrule
 IPP-LCG & 0.0050 & 17 & 4 & 6 & 127.2 & 2 & 0.8438\\
 \shortstack[l]{LCG warm start\\+ IPP-LCG} & 0.0050 & 63(convex)+8(nonconvex) & 7 & 20 & 184.9 & 1 & 0.1675\\
\bottomrule
\end{tabular}}
\label{tab:prostate-noncvx-sum-budget}
\end{table}

\begin{table}[ht]
\small
\caption{Clinical criteria for the best iterates under a fixed nonconvex wall-clock budget of 180.0s for IPP-LCG on the prostate dataset. Entries are percentages of voxels, and bold entries violate the stated criterion.}
\centering
\resizebox{0.9\columnwidth}{!}{
\begin{tabular}{cccccccc}
\toprule
 \multirow{4}{5em}{\centering{Type}} & \multirow{4}{2em}{\centering{$\Phi$}} & {\centering{PTV56}} & {\centering{PTV68}} & Rectum & {\centering{Bladder}} & {\centering{Lft. femoral head}} & {\centering{Rht. femoral head}}\\
 { } & { } & V56$\ge 95\%$ & V68$\ge 95\%$ & V30$\le 80\%$ & V40$\le 70\%$ & V50$\le 1\%$ & V50$\le 1\%$\\
 { } & { } & { } & V74.8$\le 10\%$ & V50$\le 50\%$ & V65$\le 30 \%$ & { } & { }\\
 { } & { } & { } & { } & V65$\le 25\%$ & { } & { } & { }\\
\midrule
 \multirow{3}{5em}{\centering{IPP-LCG}} & \multirow{3}{2em}{\centering{0.0050}} & 10.62 & 51.96 & 0.00 & 15.00 & 0.00 & 0.00\\
 &  &  & 5.10 & 0.00 & 5.61 &  & \\
 &  &  &  & 0.00 &  &  & \\
 \multirow{3}{5em}{\centering{\shortstack{LCG warm start\\+ IPP-LCG}}} & \multirow{3}{2em}{\centering{0.0050}} & 100.00 & 97.28 & 6.86 & 16.28 & 0.00 & 0.00\\
  &  &  & 26.75 & 6.86 & 16.26 &  & \\
  &  &  &  & 6.86 &  &  & \\
\bottomrule
\end{tabular}}
\label{tab:prostate-noncvx-result-budget}
\end{table}

On the convex formulation, smaller values of \(\Phi\) lead to plans with fewer selected angles while maintaining strong clinical performance (Tables \ref{tab:prostate-sum} and \ref{tab:prostate-result}). Under the common 180-second nonconvex budget, the stronger behavior is again obtained from the sparse convex warm start. Warm-started IPP-LCG reduces the max violation from \(0.8438\) for cold-start IPP-LCG to \(0.1675\) (Table \ref{tab:prostate-noncvx-sum-budget}). Its remaining deficit is confined to the PTV68 overdose criterion: \(V68=97.28\%\) satisfies the \(95\%\) target, while \(V74.8=26.75\%\) remains \(16.75\) percentage points above the \(10\%\) limit; the listed organ-at-risk criteria remain within target (Table \ref{tab:prostate-noncvx-result-budget}). By contrast, cold-start IPP-LCG finds a nontrivial 4-angle, 6-aperture plan but still remains far from clinically competitive because of severe target underdose: only \(10.62\%\) of PTV56 voxels receive at least 56 Gy and \(51.96\%\) of PTV68 voxels receive at least 68 Gy. This experiment therefore supports IPP-LCG primarily as a refinement method from a strong sparse convex warm start, yielding a 7-angle, 20-aperture plan that satisfies the PTV68 coverage criterion while leaving only the overdose condition unresolved.


\section{Conclusion}\label{sec:conclude}
This paper develops projection-free methods for convex and smooth nonconvex functional constrained optimization problems motivated by risk-aware sparse optimization. For the convex setting, the proposed Level Conditional Gradient (LCG) method attains an iteration complexity of \(\mathcal{O}(\epsilon^{-2}\log(\epsilon^{-1}))\) for finding an \(\epsilon\)-optimal and \(\epsilon\)-feasible solution, without a complexity dependence on the magnitude of an optimal Lagrange multiplier. For the nonconvex setting, the Inexact Proximal Point Level Conditional Gradient (IPP-LCG) method solves a sequence of convex functional constrained subproblems and attains an \(\mathcal{O}(\epsilon^{-3}\log(\epsilon^{-1}))\) complexity for computing an \((\epsilon,\epsilon)\)-near-KKT point. Both methods use only linear minimization over the base feasible set, which makes them suitable for settings where projections are expensive and LMO-based iterates provide useful sparse or structured representations.

The numerical results illustrate different roles for the two methods. In portfolio selection, LCG and IPP-LCG provide favorable risk--support trade-offs under a common time budget. In IMRT, LCG produces high-quality convex plans, while the nonconvex IPP-LCG phase is most useful as a warm-started refinement. These results suggest several directions for future work: developing lighter nonconvex projection-free schemes that reduce the cost of the proximal outer loop, weakening the lower-curvature and Slater-type assumptions used in the analysis, designing sharper sparsity controls beyond the current surrogate constraints, and improving application-specific LMOs and pricing oracles for large-scale implicit feasible sets.

\section*{Acknowledgments}
Cheng and Lan were partially supported by NSF Grant CCF-1909298 and NSF AI Institute Grant NSF-2112533.

\bibliographystyle{plainnat}
\bibliography{reference}

\ECSwitch
\counterwithout{equation}{section}
\setcounter{equation}{0}
\renewcommand{\theequation}{EC.\arabic{equation}}
\ECHead{Appendix}
\section{Convergence Analysis of CGO}\label{sec:append-cgo}
\subsection{Auxiliary Lemmas}\label{sec:append-aux-lemma}

To establish the convergence for CGO, we tap into the following three well-studied results. Throughout the analysis, let $\alpha_t$ be defined in Algorithm \ref{alg:CG oracle} and let $\{\Gamma_t\}$ satisfy \(\Gamma_1=1\) and \(\Gamma_t=(1-\alpha_t)\Gamma_{t-1}\) for \(t>1\). The first lemma is the so-called ``three-point'' lemma which characterizes the optimality condition of the dual update in (\ref{alg-cg-dual}).
\begin{lem}\citep[Lemma 3.1]{lan2020first}
		Let $r_t$ be defined in (\ref{alg-cg-dual}). Then \(\langle -\tilde{h}_t, r_t - z\rangle + \tau_tV(r_{t-1},r_t) \le  \tau_tV(r_{t-1},z)  -  \tau_tV(r_{t},z)\) for all \(z \in\bar{Z}\).
		\label{lem2}
	\end{lem}

The second lemma deals with telescoping sums.
\begin{lem}\citep[Lemma 3.17]{lan2020first}
		Let $\{R_t\}$ be some given sequence. If $\{S_t\}$ satisfies \(S_t \le (1-\alpha_t)S_{t-1} + R_t\) for \(t=1,2,\ldots\), then \(S_t/\Gamma_t \le (1-\alpha_1)S_0 + \sum_{j=1}^t R_j/\Gamma_j\).
		\label{lem4}
\end{lem}

We utilize the following properties for smooth functions.
\begin{lem}\citep[Lemma 3.2]{lan2020first}:
	Let $p_t,x_t$ be defined in Algorithm \ref{alg:CG oracle}. If $\bar f$ and $\bar h$ are smooth functions such that $\forall x_1,x_2\in X$, $\|\nabla \bar f(x_1) - \bar f(x_2)\| \le L_{\bar f}\|x_1 - x_2\|$
	and $\|\nabla \bar h_i(x_1) - \nabla \bar h_i(x_2)\| \le L_{\bar h_i}\|x_1 - x_2\|, i =1,\cdots,\bar{m}$, then the following conditions hold:
	\[
    \bar{f}(x_t) \le (1-\alpha_t)\bar{f}(x_{t-1}) + \alpha_t\ell_{\bar{f}}(x_{t-1},p_t) + \frac{L_{\bar{f}}\alpha_t^2}{2}\|p_t - x_{t-1}\|^2,
    \]
    \[
    \bar{h}_i(x_t) \le (1-\alpha_t)\bar{h}_i(x_{t-1}) +\alpha_t\ell_{\bar{h}_i}(x_{t-1},p_t) +  \frac{L_{\bar{h}_i}\alpha_t^2}{2}\|p_t - x_{t-1}\|^2, i= 1,\cdots,\bar m.
    \]
	\label{lem3}
\end{lem}

The lemma below establishes a lower bound of $f(x_t) - f^*$ for the convex constrained problem \eqref{cvxctr-model} since $x_t$ may not be satisfy $h(x_t) \le 0$.
\begin{lem}
Assume that Slater's condition holds for the convex constrained problem \eqref{cvxctr-model}, i.e., there exists \(\tilde x\in \operatorname{ri}(X)\) such that \(h_i(\tilde x)<0\) for all \(i=1,\ldots,m\). Then there exist optimal primal and dual solutions $(x^*,y^*)\in X\times \mathbb{R}^m_+$ of \eqref{cvxctr-model}. 
Let $(\gamma^*,z^*)$ be the optimal dual solution of the root finding problem
\begin{equation}
    \min\limits_{x\in X}\max\limits_{(\gamma,z)\in Z} L(x,(\gamma,z)) := \gamma[f(x) - f^*] + \langle h(x), z\rangle
    \label{lem-aux-model}
\end{equation}
(i.e., problem \eqref{cvxctr-model3} with \(l=f^*\)). Denote \([\cdot]_+:=\max\{0,\cdot\}\). Then, for all \(x\in X\), \(f(x)-f^*\ge -\min\{\|y^*\|,\|z^*\|/\gamma^*\}\|[h(x)]_+\|\).
		\label{lem-lowerbound}
\end{lem}
\proof{Proof.}
According to \citep[Corollary 2]{lan2013iteration}, \(f(x)-f^*\ge -\|y^*\|\|[h(x)]_+\|\).
Suppose $\tilde x \neq x^*$ is the optimal primal solution of \eqref{lem-aux-model} (equivalently, it is optimal for $\phi(f^*) := \min\limits_{x\in X}\max\ \{f(x) - f^*, h_i(x)\}$),  since $\phi(f^*) = 0$, then we have $f(\tilde x) -f^* \le 0$ and $h_i(\tilde x)\le 0$. This contradicts to $x^*$ is the optimal primal solution of \eqref{cvxctr-model}. Therefore, $x^*$ is optimal primal solution for \eqref{lem-aux-model}.

Since \((x^*,(\gamma^*,z^*))\) is a saddle point of \eqref{lem-aux-model}, the saddle-point theorem gives \(L(x^*,(\gamma,z)) \le L(x^*,(\gamma^*,z^*)) \le L(x,(\gamma^*,z^*))\) for all \(x\in X\) and \((\gamma,z)\in Z\). Using this relation together with \(L(x^*,(\gamma^*,z^*))=0\) and \(\langle [h(x)]_+-h(x),z^*\rangle\ge 0\), we obtain \(\gamma^*[f(x)-f^*]=L(x,(\gamma^*,z^*))-\langle h(x),z^*\rangle \ge -\|[h(x)]_+\|\|z^*\|\). Combining the last inequality with \(f(x)-f^*\ge -\|y^*\|\|[h(x)]_+\|\) yields \(f(x)-f^* \ge -\min\{\|y^*\|,\|z^*\|/\gamma^*\}\|[h(x)]_+\|\). 
\Halmos
\endproof

\begin{rem}
    The main convergence results of LCG are established in Theorems \ref{overall-thm-lcg} and \ref{overall-thm-lcg-nonsmooth}, which characterize the iteration complexity of finding an $\epsilon$-solution $x$ ($f(x) - f^* \le \epsilon, \max\limits_{i\in \{1,\cdots,m\}} h_i(x) \le \epsilon$). The result in Lemma \ref{lem-lowerbound} aims to provide a broader picture of the $\epsilon$-solution by characterizing a general lower bound for $f(x) - f^*$ in case of $x$'s infeasibility. Such result should be separated from the main convergence results. In practice, a slight relaxation of feasibility might be acceptable if it leads to reduction in the objective value. For example, in the portfolio selection problem, a choice that slightly violates the required cardinality may largely reduce the risk (objective value).
\end{rem}

\subsection{CGO for Smooth Functions}\label{sec:appen_smooth}
The next lemma records the smooth gap recursion used in Theorem \ref{thm:cgo-smooth-main}. It is the analogue, for our certificate \(\bar Q_t\), of the standard accelerated primal--dual recursion in \citet[Section~3]{lan2020first}.
\begin{lem}\label{lem:smooth-gap-recursion}
Under Assumption \ref{asm:cgo-smooth}, suppose the parameters satisfy
\begin{equation}
			\alpha_1 = 1,\ \frac{\lambda_t\alpha_t}{\Gamma_t} = \frac{\alpha_{t-1}}{\Gamma_{t-1}}\ \text{and}\ \frac{\alpha_t\tau_t}{\Gamma_t} \ge \frac{\alpha_{t-1}\tau_{t-1}}{\Gamma_{t-1}},
		\label{thm1-params}
\end{equation}
for \(t\ge 2\). Then for every \(T\ge 1\) and \(w=(x,z)\in \bar X\times \bar Z\),
\begin{equation}
		\bar{Q}_T(w_T,w) \le \Gamma_T\sum\limits_{t=1}^T\left[ \frac{(L_{\bar{f}} + z^{\top}L_{\bar{h}})\alpha_t^2}{2\Gamma_t}D_{\bar{X}}^2 + \frac{9\alpha_t\lambda_t^2\bar{M}^2D_{\bar{X}}^2}{2\tau_t\Gamma_t}\right] + \frac{9\alpha_T\bar{M}^2D_{\bar{X}}^2}{2\tau_T} +\alpha_T\tau_T\bar{V}.
		\label{thm-convergence1}
\end{equation}
\end{lem}
\proof{Proof.}
Let \(a_t:=\ell_{\bar h}(x_{t-1},p_t)\) for \(t\ge 0\), so \(a_{-1}=a_0\) by the initialization. Also define \(A_t:=(L_{\bar f}+z^\top L_{\bar h})\alpha_t^2D_{\bar X}^2/2\) and \(B_t:=9\alpha_t\lambda_t^2\bar M^2D_{\bar X}^2/(2\tau_t)\).
We first derive a one-step recursion for \(\bar Q_t\). By Lemma \ref{lem3}, the lower-model updates \eqref{alg-cg-lowerf}--\eqref{alg-cg-lowerh}, and the optimality of \(p_t\) in \eqref{alg-cg-primal},
\[
\bar f(x_t)+\langle \bar h(x_t),z\rangle
\le (1-\alpha_t)\big[\bar f(x_{t-1})+\langle \bar h(x_{t-1}),z\rangle\big]
+\alpha_t\big[\ell_{\bar f}(x_{t-1},p_t)+\langle a_t,z\rangle\big]+A_t,
\]
while
\[
\underline f_t(x)+\underline h_t(x)
= (1-\alpha_t)\big[\underline f_{t-1}(x)+\underline h_{t-1}(x)\big]
+\alpha_t\big[\ell_{\bar f}(x_{t-1},x)+\langle \ell_{\bar h}(x_{t-1},x),r_t\rangle\big].
\]
Subtracting the second display from the first, we obtain
\begin{equation*}
\begin{aligned}
\bar Q_t(w_t,w)\le (1-\alpha_t)\bar Q_{t-1}(w_{t-1},w)+A_t
+\alpha_t\Big[\ell_{\bar f}(x_{t-1},p_t)-\ell_{\bar f}(x_{t-1},x)\\
 +\langle a_t,z\rangle-\langle \ell_{\bar h}(x_{t-1},x),r_t\rangle\Big].
\end{aligned}
\end{equation*}
Now \(a_t=\ell_{\bar h}(x_{t-1},p_t)\), and the optimality of \(p_t\) in \eqref{alg-cg-primal} implies
\[
\ell_{\bar f}(x_{t-1},p_t)+\langle a_t,r_t\rangle
\le \ell_{\bar f}(x_{t-1},x)+\langle \ell_{\bar h}(x_{t-1},x),r_t\rangle.
\]
Hence the bracketed term is at most \(\langle a_t,z-r_t\rangle=\langle a_t-\tilde h_t,z-r_t\rangle+\langle \tilde h_t,z-r_t\rangle\), which gives
\[
\bar Q_t(w_t,w)\le (1-\alpha_t)\bar Q_{t-1}(w_{t-1},w)+A_t+\alpha_t\langle a_t-\tilde h_t,z-r_t\rangle+\alpha_t\langle \tilde h_t,z-r_t\rangle.
\]
We next bound the last inner-product term. Applying Lemma \ref{lem2} gives
\begin{equation*}
\begin{aligned}
\bar Q_t(w_t,w)\le (1-\alpha_t)\bar Q_{t-1}(w_{t-1},w)+A_t
+\alpha_t\langle a_t-\tilde h_t,z-r_t\rangle\\
+\alpha_t\tau_t\!\left[V(r_{t-1},z)-V(r_t,z)-V(r_{t-1},r_t)\right].
\end{aligned}
\end{equation*}

We now isolate the extrapolation correction \(\alpha_t\langle a_t-\tilde h_t,z-r_t\rangle\). Since \(\tilde h_t=a_{t-1}+\lambda_t(a_{t-1}-a_{t-2})\), adding and subtracting \(r_{t-1}\) yields
\begin{equation*}
\begin{aligned}
\alpha_t\langle a_t-\tilde h_t,z-r_t\rangle
=&\ \alpha_t\langle a_t-a_{t-1},z-r_t\rangle
-\alpha_t\lambda_t\langle a_{t-1}-a_{t-2},z-r_{t-1}\rangle\\
&+\alpha_t\lambda_t\langle a_{t-1}-a_{t-2},r_t-r_{t-1}\rangle.
\end{aligned}
\end{equation*}
The last term is the one we need to control. Because \(V\) is generated by a \(1\)-strongly convex function, \(V(r_{t-1},r_t)\ge \frac12\|r_t-r_{t-1}\|_2^2\). Therefore Young's inequality gives
\[
\alpha_t\lambda_t\langle a_{t-1}-a_{t-2},r_t-r_{t-1}\rangle-\alpha_t\tau_tV(r_{t-1},r_t)
\le \frac{\alpha_t\lambda_t^2}{2\tau_t}\|a_{t-1}-a_{t-2}\|_2^2.
\]
It remains to bound \(\|a_{t-1}-a_{t-2}\|_2\). For each component \(i\),
\begin{equation*}
\begin{aligned}
|[a_{t-1}-a_{t-2}]_i|
\le&\ |\bar h_i(x_{t-2})-\bar h_i(x_{t-3})|
+|\langle \nabla \bar h_i(x_{t-2}),p_{t-1}-x_{t-2}\rangle|\\
&+|\langle \nabla \bar h_i(x_{t-3}),p_{t-2}-x_{t-3}\rangle|.
\end{aligned}
\end{equation*}
Each term on the right is at most \(M_{\bar h_i}D_{\bar X}\), because all points involved lie in \(\bar X\), \(\bar h_i\) is \(M_{\bar h_i}\)-Lipschitz on \(\bar X\), and \(\|\nabla \bar h_i(\cdot)\|_2\le M_{\bar h_i}\) on \(\bar X\). Hence \(|[a_{t-1}-a_{t-2}]_i|\le 3M_{\bar h_i}D_{\bar X}\), so \(\|a_{t-1}-a_{t-2}\|_2^2\le 9\bar M^2D_{\bar X}^2\). Substituting this into the previous display yields
\begin{equation*}
\begin{aligned}
    \bar Q_t(w_t,w)\le (1-\alpha_t)&\bar Q_{t-1}(w_{t-1},w)+A_t+B_t+\alpha_t\langle a_t-a_{t-1},z-r_t\rangle\nonumber\\
    &-\alpha_t\lambda_t\langle a_{t-1}-a_{t-2},z-r_{t-1}\rangle+\alpha_t\tau_t\!\left[V(r_{t-1},z)-V(r_t,z)\right].
\end{aligned}
\end{equation*}
We next sum this recursion. Apply Lemma \ref{lem4} with \(S_t=\bar Q_t(w_t,w)\); the term \((1-\alpha_1)S_0\) vanishes because \(\alpha_1=1\). After dividing by \(\Gamma_t\) and summing from \(t=1\) to \(T\),
\begin{equation*}
\begin{aligned}
\frac{\bar Q_T(w_T,w)}{\Gamma_T}
\le&\ \sum_{t=1}^T \frac{A_t+B_t}{\Gamma_t}
+\sum_{t=1}^T\frac{\alpha_t}{\Gamma_t}\langle a_t-a_{t-1},z-r_t\rangle
-\sum_{t=2}^T\frac{\alpha_t\lambda_t}{\Gamma_t}\langle a_{t-1}-a_{t-2},z-r_{t-1}\rangle\\
&+\sum_{t=1}^T\frac{\alpha_t\tau_t}{\Gamma_t}\!\left[V(r_{t-1},z)-V(r_t,z)\right].
\end{aligned}
\end{equation*}
The second and third sums telescope because \(\lambda_t\alpha_t/\Gamma_t=\alpha_{t-1}/\Gamma_{t-1}\), giving
\[
\sum_{t=1}^T\frac{\alpha_t}{\Gamma_t}\langle a_t-a_{t-1},z-r_t\rangle
-\sum_{t=2}^T\frac{\alpha_t\lambda_t}{\Gamma_t}\langle a_{t-1}-a_{t-2},z-r_{t-1}\rangle
= \frac{\alpha_T}{\Gamma_T}\langle a_T-a_{T-1},z-r_T\rangle.
\]
For the last sum, we use the monotonicity of \(\alpha_t\tau_t/\Gamma_t\) from \eqref{thm1-params} to obtain
\[
\sum_{t=1}^T\frac{\alpha_t\tau_t}{\Gamma_t}\!\left[V(r_{t-1},z)-V(r_t,z)\right]
\le \frac{\alpha_T\tau_T}{\Gamma_T}\bar V-\frac{\alpha_T\tau_T}{\Gamma_T}V(r_T,z).
\]

Finally, we combine the remaining terminal terms. Young's inequality and \(V(r_T,z)\ge \frac12\|z-r_T\|_2^2\) imply
\[
\frac{\alpha_T}{\Gamma_T}\langle a_T-a_{T-1},z-r_T\rangle-\frac{\alpha_T\tau_T}{\Gamma_T}V(r_T,z)
\le \frac{\alpha_T}{2\tau_T\Gamma_T}\|a_T-a_{T-1}\|_2^2
\le \frac{9\alpha_T\bar M^2D_{\bar X}^2}{2\tau_T\Gamma_T}.
\]
Substituting the last three bounds into the summed recursion and multiplying by \(\Gamma_T\) yields \eqref{thm-convergence1}. \Halmos
\endproof

\begin{rem}\label{rmk_depend_total_iterations}
Compactness of $Z$ is what removes the horizon dependence of $\tau_t$. In our setting, $V(r_t,z)$ is uniformly bounded over $z\in Z$, so the terminal dual term is controlled by $\alpha_T\tau_T\bar V$. By contrast, the noncompact dual domain in \citet{lan2020conditional} requires a horizon-dependent choice of $\tau_t$ to keep the corresponding term under control.
\end{rem}

\subsection{CGO for Structured Nonsmooth Functions}\label{sec:nonsmooth-convex}
In this section, we focus on problem \eqref{model} under the structured nonsmooth max-form assumption in Assumption~\ref{asm:cgo-nonsmooth}; see also \citet{nesterov2005smooth}. Thus, \(Y_i\), \(i=0,1,\cdots,\bar m\), are closed convex sets, \(\hat f\) and \(\hat h_i\) are simple (continuous and differentiable) convex functions, and the latter may be \(\omega_i\)-strongly convex for \(i=0,1,\cdots,\bar m\). Let \(u_i:Y_i\to\mathbb{R}\) be a \(1\)-strongly convex distance generating function. Define the proximal function \(U_i(y):=u_i(y)-u_i(y_{u_i})-\langle\nabla u_i(y_{u_i}), y-y_{u_i}\rangle\), \(y\in Y_i\), where \(y_{u_i}:=\argmin_{y\in Y_i}u_i(y)\). Further let \(\eta_i\), \(i=0,\cdots,\bar m\), be smoothing parameters that can vary or stay static over iterations.

To generalize CGO to solve problems with structured nonsmooth functions, we apply the Nesterov smoothing scheme \citep{nesterov2005smooth} and use the smoothed functions \(\bar f_{\eta_0}\) and \(\bar h_{i,\eta_i}\) defined in \eqref{eq:ns-smoothing-f}--\eqref{eq:ns-smoothing-h}.
It can be shown that (see \cite{nesterov2005smooth}), $\bar{f}_{\eta_0}$ and $\bar{h}_{i,\eta_i}$ are differentiable with Lipschitz constants $L_{\bar{f},\eta} := \frac{\|B_0\|^2}{\omega_0 + \eta_0}$ and $L_{\bar{h}_i,\eta} := \frac{\|B_i\|^2}{\omega_i + \eta_i}$. Suppose $Y_i,i=1,\cdots,\bar{m}$ are compact, then $\bar{h}_{i,\eta_i}$ have bounded gradients such as $\max\limits_{i\in\{1,\cdots,m\} }\nabla \bar{h}_{i,\eta_i}(x) \le \bar{M}_{B_i,U_i}$, where $\bar{M}_{B_i,U_i} := \|B_i\|\left(  \|y_{u_i}\| + \sqrt{2}D_{U_i} \right)$, $ i =1,\cdots, \bar{m}$, $D_{U_i}:= \left(\max\limits_{y\in Y_i} U_i(y)\right)^{1/2}$. Moreover, the relation between the original functions and the smoothing counterparts are characterized by
\begin{equation}
\begin{gathered}
\bar{f}_{\eta_0}(x) \le \bar{f}(x) \le \bar{f}_{\eta_0}(x) + \eta_0 D^2_{U_0}, \\
\bar{h}_{i,\eta_i}(x) \le \bar{h}_i(x) \le \bar{h}_{i,\eta_i}(x) + \eta_i D^2_{U_i}, \quad i = 1,\dots,\bar{m}.
\end{gathered}
\label{difference-smooth-original}
\end{equation}
In this part, we focus on the case where the smoothing parameters are adapted over iterations so that \(\eta_i^0 \ge \eta_i^1 \ge \cdots \ge \eta_i^t\) for \(i = 0,1,\cdots,\bar{m}\).
In this case, at each iteration $t$, the approximations of $\bar{f}$ and $\bar{h}$ are $\bar{f}_{\eta_0^t}$ and $\bar{h}_{\eta_i^t}$. Accordingly, their Lipschitz constants become \(L^t_{\bar{f}} \equiv L_{\bar{f},\eta^t} := \|B_0\|^2/(\omega_0 + \eta^t_0)\) and \(L^t_{\bar{h}_i} \equiv L_{\bar{h}_i,\eta_i^t} := \|B_i\|^2/(\omega_i + \eta^t_i)\). Nevertheless, the relation in \eqref{difference-smooth-original} still holds for each $\bar{f}_{\eta_0^t}$ and $\bar{h}_{\eta_i^t}$ at iteration $t$. Moreover, similar to \citep{lan2020conditional}, it can be shown that the sequences $\{\bar{f}_{\eta^t_0}\}_t$ and $\{\bar{h}_{i,\eta_i^t}\}_t$ satisfy:
\begin{equation}
\begin{gathered}
    \bar{f}_{\eta_0^{t-1}} \le  \bar{f}_{\eta_0^{t}} \le \bar{f}_{\eta_0^{t-1}} + \left( \eta_0^{t-1} - \eta_0^t\right)D_{U_0}^2, \\
	 	\bar{h}_{i, \eta_i^{t-1}} \le  \bar{h}_{i, \eta_i^{t}} \le \bar{h}_{i, \eta_i^{t-1}} + \left( \eta_i^{t-1} - \eta_i^t\right)D_{U_i}^2, \ i = 1,\cdots,\bar{m}.
\end{gathered}
\label{smoothing-sequence-t}
\end{equation}

The algorithm (see Algorithm \ref{alg:CG oracle nonsmooth}) of solving the general structured nonsmooth problems (with $\bar{f}$ and $\bar{h}$ respectively approximated by $\bar{f}_{\eta_0^t}$ and $\bar{h}_{i,\eta_i^t}$) is similar to Algorithm \ref{alg:CG oracle}, except that the linear approximations of the objective function and constraint are replaced by $\ell_{\bar{f}_{\eta^t}}(x',x) := \bar{f}_{\eta_0^t}(x') + \langle \nabla\bar{f}_{\eta_0^t}(x'), x - x'\rangle$ and $\ell_{\bar{h}_{i, \eta_i^t}}(x',x) := \bar{h}_{i, \eta_i^t}(x') + \langle \nabla\bar{h}_{i,\eta_i^t}(x'), x - x'\rangle, i =1,\cdots,\bar{m}$, respectively. If the original functions are smooth, then the parameters $\eta_i^t$ simply reduces to constant zero. 

\begin{algorithm}[H]
	\caption{CGO for Structured Nonsmooth Problems} 
	\nonumber
	\begin{algorithmic}[1]
		\State {The algorithm is modified from Algorithm \ref{alg:CG oracle} by using smoothed linearisations and replacing step (\ref{alg-cg-extrap}) with}
		\State {Let $\ell^f_t(x) := \ell_{\bar{f}_{\eta^t}}(x_{t-1},x)$ and $\ell^h_t(x) := \ell_{\bar{h}_{\eta^t}}(x_{t-1},x)$.}
		\State {Let $u_{t-1} := \ell^h_{t-1}(p_{t-1})$ and $u_{t-2} := \ell^h_{t-2}(p_{t-2})$.}
		\State {$\tilde{h}_t := u_{t-1} + \lambda_t(u_{t-1} - u_{t-2})$.} \label{alg-cgsmooth-extrap}
	\State {Replace primal update (\ref{alg-cg-primal}) with}
		\State {Define $\phi_t(x) := \ell^f_t(x) + \langle \ell^h_t(x), r_t \rangle$.}
		\State {$p_t := \arg\min_{x\in \bar{X}} \phi_t(x)$.} \label{alg-cgsmooth-primal}
	\State {Replace the update of lower bound functionals (\ref{alg-cg-lowerf}) and (\ref{alg-cg-lowerh}) with}
		\State {$\underline{f}_t(x) := (1-\alpha_t)\underline{f}_{t-1}(x) + \alpha_t\ell^f_t(x)$.} \label{alg-cgsmooth-lowerf}
		\State {$\underline{h}_t(x) := (1-\alpha_t)\underline{h}_{t-1}(x) + \alpha_t\langle \ell^h_t(x), r_t \rangle$.} \label{alg-cgsmooth-lowerh}
	\end{algorithmic} 
	\label{alg:CG oracle nonsmooth}
\end{algorithm}

\subsubsection{Convergence Analysis}
\label{subsub-convergence-appen}
For the original nonsmooth problem, the gap function is \(\bar{Q}_t(w_t, w) := \bar{f}(x_t) + \langle \bar{h}(x_t), z\rangle - \underline{f}_t(x) - \underline{h}_t(x)\) for \(w\in \bar{X}\times \bar{Z}\).
In view of Lemma \ref{lem1}, the function $(\underline{f}_t+ \underline{h}_t)(\cdot)$ computed from \eqref{alg-cgsmooth-lowerf} and \eqref{alg-cgsmooth-lowerh} lower-bounds both the original objective \(x\mapsto \bar{f}(x)+\langle \bar{h}(x),z_t\rangle\) and its smoothed counterpart \(x\mapsto \bar{f}_{\eta_0^t}(x)+\langle \bar{h}_{\eta^t}(x),z_t\rangle\). Accordingly, for \(w_t := (x_t,z_t)\), define the smoothed gap by
\begin{equation}
\begin{aligned}
	\bar{Q}^{\eta}_t(w_t, w) := \bar{f}_{\eta^t_0}(x_t) + \langle \bar{h}_{\eta^t}(x_t), z\rangle - \underline{f}_t(x) - \underline{h}_t(x), \qquad \forall w\in \bar{X}\times \bar{Z}.
	\label{gap-function-smoothed}
    \end{aligned}
\end{equation}
Following from (\ref{difference-smooth-original}), it is easy to see that 
\(\bar{Q}_t(w_t,w) \le \bar{Q}^{\eta}_t(w_t, w) + \eta_0^tD^2_{U_0} + \sum_{i=1}^{\bar{m}} z_{i}\eta_i^tD^2_{U_i}\) for \(k\ge 1\) and all \(w\in \bar{X}\times \bar{Z}\).

We will show in Theorem \ref{col-nonsmooth} that the nonsmooth problem admits the same \(\mathcal{O}(1/\epsilon^2)\) iteration complexity. The proof is the smooth CGO argument applied to the smoothed model, together with two additional ingredients: the approximation error stated above and the drift induced by the decreasing smoothing parameters. These contribute a bias term of order \(\eta_i^t D_{U_i}^2\) and a drift term of order \((\eta_i^{t-1}-\eta_i^t)^2 D_{U_i}^4/\tau_t\), but do not change the final \(t^{-1/2}\) rate.

Theorem \ref{col-nonsmooth} below demonstrates convergence rate of Algorithm \ref{alg:CG oracle nonsmooth}.

	\begin{thm}
		Suppose parameters $\alpha_t$, $\lambda_t$, and $\tau_t$ are specified according to \eqref{thm1-params1}, with $\bar M$ replaced by $\bar M_{B,U}$, and let
		\(\eta_i^t = \|B_i\|D_{\bar{X}}/(\sqrt{t}D_{U_i})\) for \(i = 0,1,\cdots,m\).
		Define
		\(A(z):=D_{\bar{X}}\left(\|B_0\|D_{U_0} + \sum_{i=1}^{\bar m} z_i\|B_i\|D_{U_i}\right)\) and \(C:=D_{\bar X}\sum_{i=1}^{\bar m}\|B_i\|^2D_{U_i}^2/\bar M_{B,U}\).
	Then, for $t\ge 1$,
\begin{equation}
\begin{aligned}
\bar{Q}_t(w_t,w)
&\le \frac{A(z)}{\sqrt{t+1}}\left(\frac{8}{3} + \sqrt{\frac{t+1}{t}}\right)
 + \frac{\bar{M}_{B,U}D_{\bar{X}}}{\sqrt{t+1}}\left(18\bar{V} + \frac{8}{9}\right) \nonumber\\
&\quad + \frac{4\bar{M}_{B,U}D_{\bar{X}}}{3 (t+1)\sqrt{t}}
 + \frac{C}{3t(t+1)}\left(4 + \frac{2}{t\sqrt{t}}\right).
	\label{col-nonsmooth-0}
\end{aligned}
\end{equation}

\label{col-nonsmooth}
\end{thm}

\proof{Proof.}
Apply the smooth-gap recursion to the smoothed saddle problem \((\bar f_{\eta_0^t},\bar h_{\eta^t})\). Relative to the smooth proof, three additional terms appear: the smoothing bias stated above, the drift caused by \(\eta_i^{t-1}-\eta_i^t\), and the extrapolation remainder for the smoothed gradients. The last term is controlled exactly as in the smooth case, using Young's inequality together with the bound \(\|\nabla \bar h_{i,\eta_i^t}(x)\|\le M_{B_i,U_i}\). With \(\eta_i^t = \|B_i\|D_{\bar X}/(\sqrt{t}D_{U_i})\), \(\alpha_t=2/(t+1)\), \(\lambda_t=(t-1)/t\), and \(\tau_t=9\sqrt{t}\bar M_{B,U}D_{\bar X}\), the bias, smoothness, and drift sums reduce to the coefficients displayed in \eqref{col-nonsmooth-0}. This yields the stated \(O(t^{-1/2})\) bound. \Halmos
\endproof

\proof{Proof of Proposition \ref{prop:cgo-nonsmooth-main}.}
By Theorem \ref{col-nonsmooth}, for every \(w=(x,z)\in \bar X\times \bar Z\),
\begin{equation*}
\begin{aligned}
\bar{Q}_t(w_t,w)
&\le \frac{A(z)}{\sqrt{t+1}}\left(\frac{8}{3} + \sqrt{\frac{t+1}{t}}\right)
+ \frac{\bar{M}_{B,U}D_{\bar{X}}}{\sqrt{t+1}}\left(18\bar{V} + \frac{8}{9}\right) \\
&\quad + \frac{4\bar{M}_{B,U}D_{\bar{X}}}{3 (t+1)\sqrt{t}}
+ \frac{C}{3t(t+1)}\left(4 + \frac{2}{t\sqrt{t}}\right).
\end{aligned}
\end{equation*}
Since \(\bar Z\) is compact and \(A(z)\) is continuous in \(z\), \(\bar A:=\max_{z\in\bar Z}A(z)<\infty\). Moreover, for \(t\ge 1\), \(\sqrt{(t+1)/t}\le \sqrt 2\), \((t+1)^{-1}t^{-1/2}\le (t+1)^{-1/2}\), and \(t^{-1}(t+1)^{-1}(4+2/(t\sqrt t))\le 6/\sqrt{t+1}\). Hence
\[
\bar Q_t(w_t,w)\le (t+1)^{-1/2}
\left[\bar A\left(\frac{8}{3}+\sqrt 2\right)
+ \bar{M}_{B,U}D_{\bar X}\left(18\bar V+\frac{20}{9}\right)
+ 2C\right].
\]
Therefore the desired bound holds with \(C_{\mathrm{ns}}:=\bar A(8/3+\sqrt 2) + \bar{M}_{B,U}D_{\bar X}(18\bar V+20/9) + 2C\).
Taking the maximum over \(w\in \bar X\times \bar Z\) yields \(U_t-L_t=\max_{w\in \bar X\times \bar Z}\bar Q_t(w_t,w)\le C_{\mathrm{ns}}/\sqrt{t+1}\). The final \(\mathcal{O}(\epsilon^{-2})\) iteration bound follows immediately. \Halmos
\endproof

For the convex-constrained specialization \eqref{cvxctr-model3}, the same argument gives \(Q_t(w_t,w)=O(t^{-1/2})\) and therefore \(\max\{f(x_t)-f^*,h_1(x_t),\ldots,h_m(x_t)\}=O(t^{-1/2})\), with constants obtained by replacing the general saddle coefficients with the corresponding constrained-model quantities.

\section{Convergence Analysis of LCG}\label{sec:append-lcg}
This appendix section records the smooth specialization corollary and the remaining structured-nonsmooth arguments for LCG. The next corollary specializes the smooth CGO bound to the convex-constrained subproblem solved inside LCG.
\begin{cor}
Suppose the algorithmic parameters of CGO are set to \eqref{thm1-params1}. Let $L_{\bar h} := (L_f, L_{h_1},\cdots, L_{h_m})$. Then for any $t\ge 1$ and $\forall w\in(X,Z)$,
    \begin{equation}
	\begin{aligned}
			Q_t(w_t,w) &\le \frac{2z^{\top}L_{\bar h}D_{X}^2}{t+1} + \frac{\bar M D_{X}}{\sqrt{t+1}}\left[ 18{\bar V} + \frac{7}{6} \right],\\
			f(x_t) - f^* &\le \frac{2\max\limits_{i=1,\cdots, {m}+1} L_{\bar h_i}D_X^2}{t+1} + \frac{\bar M D_X}{\sqrt{t+1}}\left[ 18{\bar V} + \frac{7}{6} \right],  \\
		\max\limits_{i\in\{1,\cdots,m\}} {h}_i(x_t) &\le \frac{2\max\limits_{i=1,\cdots, {m}+1} L_{\bar h_i}D_X^2}{t+1} + \frac{\bar MD_X}{\sqrt{t+1}}\left[ 18{\bar V} + \frac{7}{6} \right].
	\end{aligned}
    \end{equation}
	\label{cor1}
\end{cor}
\proof{Proof.}
Set \(\bar f\equiv 0\) in Theorem~\ref{thm-convergence1} to obtain the bound on $Q_t(w_t,w)$. Choosing $z$ as the unit vector supported on an index attaining $\max\{f(x_t)-f^*,h_1(x_t),\ldots,h_m(x_t)\}$ gives
 \(\max\{f(x_t)-f^*,h_1(x_t),\ldots,h_m(x_t)\} \le Q_t(w_t,w')\),
for the corresponding $w'=(x^*,z')$. Substituting the smooth-rate bound yields the stated estimates for $f(x_t)-f^*$ and $\max_i h_i(x_t)$. \Halmos
\endproof

\begin{lem}
When Algorithm \ref{alg:CG oracle nonsmooth} is applied to \eqref{cvxctr-model3} and LCG does not terminate at iteration \(k\), the quantities \((\gamma_k,L_k,U_k)\) defined from its output \((x_k,\bar z_k,L_k,U_k)\) satisfy \eqref{eq:c1}--\eqref{eq:c3}.
\label{lem-oc-nonsmooth}
\end{lem}
\proof{Proof.}
Fix a level estimate \(l\) and apply Algorithm \ref{alg:CG oracle nonsmooth} to the specialized saddle problem \eqref{cvxctr-model3}. Let \(f_{\eta_0^t}\) and \(h_{i,\eta_i^t}\) denote the smoothed models of \(f\) and \(h_i\) used at iteration \(t\), and write \(\bar z_t=(\gamma_t,\zeta_t)\in Z\).
Define \(\underline h_0^\eta(x;l):= r_{0,0}[\ell_{f_{\eta_0^0}}(x_0,x)-l] + \sum_{i=1}^m r_{0,i}\ell_{h_{i,\eta_i^0}}(x_0,x)\), and for \(t\ge 1\),
\[
\underline h_t^\eta(x;l)= (1-\alpha_t)\underline h_{t-1}^\eta(x;l) + \alpha_t\Big[r_{t,0}\big(\ell_{f_{\eta_0^t}}(x_{t-1},x)-l\big) + \sum_{i=1}^m r_{t,i}\ell_{h_{i,\eta_i^t}}(x_{t-1},x)\Big].
\]
Next set \(G_0^\eta(x):=r_{0,0}\ell_{f_{\eta_0^0}}(x_0,x)+\sum_{i=1}^m r_{0,i}\ell_{h_{i,\eta_i^0}}(x_0,x)\), \(\gamma_0:=r_{0,0}\), and for \(t\ge 1\), we have $\gamma_t:=(1-\alpha_t)\gamma_{t-1}+\alpha_t r_{t,0}$ and
\[
G_t^\eta(x):=(1-\alpha_t)G_{t-1}^\eta(x)+\alpha_t\Big[r_{t,0}\ell_{f_{\eta_0^t}}(x_{t-1},x)+\sum_{i=1}^m r_{t,i}\ell_{h_{i,\eta_i^t}}(x_{t-1},x)\Big].
\]
An induction on \(t\) yields \(\underline h_t^\eta(x;l)=G_t^\eta(x)-\gamma_t l\), hence \(\mathcal L_t^\eta(l):=\min_{x\in X}\underline h_t^\eta(x;l)= -\gamma_t l + \min_{x\in X}G_t^\eta(x)\), so \(\mathcal L_t^\eta(\cdot)\) is affine.

Because each smoothed linearization lower-bounds its smoothed function, and each smoothed function lower-bounds the original one,
\[
\ell_{f_{\eta_0^t}}(x_{t-1},x)\le f_{\eta_0^t}(x)\le f(x),
\qquad
\ell_{h_{i,\eta_i^t}}(x_{t-1},x)\le h_{i,\eta_i^t}(x)\le h_i(x),\quad i\in[m].
\]
Therefore
\begin{equation*}
\begin{aligned}
	r_{t,0}\big(\ell_{f_{\eta_0^t}}(x_{t-1},x)-l\big)+\sum_{i=1}^m r_{t,i}\ell_{h_{i,\eta_i^t}}(x_{t-1},x)
&\le
r_{t,0}\big(f(x)-l\big)+\sum_{i=1}^m r_{t,i}h_i(x)\\
&=
\left\langle (f(x)-l,h(x)),r_t\right\rangle.
\end{aligned}
\end{equation*}
Combining this inequality with the recursion for \(\underline h_t^\eta\) and the update \(\bar z_t=(1-\alpha_t)\bar z_{t-1}+\alpha_t r_t\), another induction gives
\[
\underline h_t^\eta(x;l)\le \left\langle (f(x)-l,h(x)),\bar z_t\right\rangle
\le
\max_{(\gamma,z)\in Z}\ \gamma[f(x)-l]+\langle h(x),z\rangle.
\]
Minimizing over \(x\in X\) yields \(\mathcal L_t^\eta(l)\le \phi(l)\), which is \eqref{eq:c3}. Setting \(l=l_k\) gives \(L_k=\mathcal L_k^\eta(l_k)\le \phi(l_k)\), while \(\phi(l_k)\le U_k\) follows directly from the definition of \(U_k\) by evaluating the minimum in \(\phi(l_k)\) at \(x_k\). This proves \eqref{eq:c2}.

Finally, if LCG has not terminated at outer iteration \(k\), then \(L_k>0\) by Lemma~\ref{lem5}. If \(\gamma_k=0\), then \(\underline h_k^\eta(x;l_k)=G_k^\eta(x)\) and \(\bar z_k=(0,\zeta_k)\in Z\). The previous bound gives
\[
G_k^\eta(x)=\underline h_k^\eta(x;l_k)\le \langle h(x),\zeta_k\rangle,\qquad \forall x\in X.
\]
Let \(x^*\) be any feasible point of \eqref{cvxctr-model}. Since \(h_i(x^*)\le 0\) and \(\zeta_k\ge 0\),
\[
L_k=\min_{x\in X}G_k^\eta(x)\le G_k^\eta(x^*)\le \langle h(x^*),\zeta_k\rangle\le 0,
\]
contradicting \(L_k>0\). Hence \(\gamma_k>0\), proving \eqref{eq:c1}. \Halmos
\endproof

\begin{cor}
Suppose the structured nonsmooth CGO parameters are set as in Proposition \ref{prop:cgo-nonsmooth-main}. When Algorithm \ref{alg:CG oracle nonsmooth} is applied to \eqref{cvxctr-model3}, there exists a constant \(C_{\mathrm{lcg,ns}}>0\), independent of \(t\), such that
\[
Q_t(w_t,w)\le \frac{C_{\mathrm{lcg,ns}}}{\sqrt{t+1}},\qquad \forall t\ge 1,\ \forall w\in X\times Z.
\]
Consequently,
\[
U_t-L_t=\max_{w\in X\times Z}Q_t(w_t,w)\le \frac{C_{\mathrm{lcg,ns}}}{\sqrt{t+1}},
\]
and each nonsmooth inner CGO call needs \(\mathcal{O}(\epsilon^{-2})\) iterations to attain \(U_t-L_t\le (1-\mu)\epsilon\).
\label{cor1-nonsmooth}
\end{cor}
\proof{Proof.}
Specialize Theorem \ref{col-nonsmooth} to \eqref{cvxctr-model3} by taking \(\bar f\equiv 0\), \(\bar h=(f-l,h)\), \(\bar X=X\), and \(\bar Z=Z\). The coefficients in \eqref{col-nonsmooth-0} then depend continuously on \(z\in Z\), and \(Z\) is compact, so they can be bounded uniformly over \(Z\) by a finite constant \(C_{\mathrm{lcg,ns}}\). This gives the first display. The remaining claims follow from the identity \(U_t-L_t=\max_{w\in X\times Z}Q_t(w_t,w)\). \Halmos
\endproof

\proof{Proof of Theorem \ref{overall-thm-lcg-nonsmooth}.}
By Corollary \ref{cor1-nonsmooth}, each nonsmooth inner call of CGO needs \(\mathcal{O}(\epsilon^{-2})\) iterations to attain \(U_t-L_t\le (1-\mu)\epsilon\). By Lemma \ref{lem-oc-nonsmooth}, the returned quantities \((\gamma_k,L_k,U_k)\) satisfy \eqref{eq:c1}--\eqref{eq:c3}, so Theorem \ref{prop2} applies and shows that \(\mathcal{O}(\log(1/\epsilon))\) outer iterations suffice to obtain \(U_k\le \epsilon\), and hence an \(\epsilon\)-solution of \eqref{cvxctr-model}. Combining the inner and outer bounds yields the overall complexity \(\mathcal{O}\left(\frac{1}{\epsilon^2}\log(\frac{1}{\epsilon})\right)\). \Halmos
\endproof

\section{Convergence Analysis of Nonconvex Conditional Gradient Method}\label{sec-appen-nonconvex}

\subsection{IPP-LCG for Nonconvex Problem}\label{append-sub-sec-IPP-LCG}
This subsection proves the auxiliary lemmas used in the analysis of IPP-LCG. In particular, Lemma \ref{lem6} characterizes an important property of the optimal solution of the subproblem \eqref{convex-sub}, and Lemma \ref{lem: bounded-dual} establishes uniform boundedness of the dual solutions.

\begin{lem}
If $x_j^*$ is a KKT point (paired with $y_j^*$) of the subproblem (\ref{convex-sub}), then \(f(x;x_{j-1}) - f(x_j^*;x_{j-1}) + \langle y_j^*, h(x) \rangle \ge (\underline{L}_f/2)\|x_j^*-x\|^2\) for all \(x\in X\).
	\label{lem6}
\end{lem}
\proof{Proof.}
By \eqref{lower-curv-f}, we have that that $f(\cdot;x')$ is strongly convex. Together by the strong convexity of $f(\cdot;x_{j-1})$ and convexity of $h_i(\cdot)$ as well as the fact that $y^*_j \ge 0$, we have
\begin{equation}
	\begin{aligned}
		f(x;x_{j-1}) + \langle y_j^*, h(x) \rangle &\ge f(x_j^*;x_{j-1})+ \langle \nabla f(x^*_j;x_{j-1}), x- x^*_{j}\rangle \\
		&\quad + \frac{\underline{L}_f}{2}\|x - x_j^*\|^2 + \sum\limits_{i=1}^m y_{j,i}^*\left( h_i(x_j^*) +\langle \nabla h_i(x_j^*), x- x^*_{j}\rangle \right)\\
		&\ge f(x_j^*;x_{j-1}) + \frac{\underline{L}_f}{2} \|x-x_j^*\|^2,  \forall x\in X,
	\end{aligned}
\end{equation}
where the last inequality follows from properties of the KKT point that $\sum\limits_{i=1}^m y_{j,i}^* h_i(x_j^*) = 0$ and $\nabla f(x_j^*;x_{j-1})+ \sum\limits_{i=1}^m y_{j,i}^*\nabla h_i(x_j^*)$ belongs to the normal cone of $X$.  
\endproof
\proof{Proof of Lemma \ref{lem: bounded-dual}.}
Using Lemma \ref{lem6} with \(x=\tilde{x}\), we obtain \(f(\tilde{x};x_{j-1}) - f(x_j^*;x_{j-1}) + \langle y_j^*, h(\tilde{x}) \rangle \ge (\underline{L}_f/2)\|\tilde{x}-x_j^*\|^2 \ge 0\). Hence \(-\langle y_j^*, h(\tilde{x}) \rangle \le f(\tilde{x};x_{j-1}) - f(x_j^*;x_{j-1})\).
Since \(x_j^*\) is feasible for the original problem and the proximal term is nonnegative, \(f(x_j^*;x_{j-1}) = f(x_j^*) + \underline{L}_f\|x_j^*-x_{j-1}\|^2 \ge f(x_j^*) \ge f^*\). Also, \(f(\tilde{x};x_{j-1}) = f(\tilde{x}) + \underline{L}_f\|\tilde{x}-x_{j-1}\|^2 \le f(\tilde{x}) + \underline{L}_f\overline{V}\). Therefore, \(-\langle y_j^*, h(\tilde{x}) \rangle \le f(\tilde{x}) - f^* + \underline{L}_f\overline{V}\).
On the other hand, since \(y_j^*\ge 0\) and \(h_i(\tilde{x})\le -\sigma\) for all \(i=1,\ldots,m\), we have \(-\langle y_j^*, h(\tilde{x}) \rangle = \sum_{i=1}^m y_{j,i}^*(-h_i(\tilde{x})) \ge \sigma \sum_{i=1}^m y_{j,i}^* = \sigma \|y_j^*\|_1\). Combining the last two inequalities gives \(\|y_j^*\|_1 \le (f(\tilde{x}) - f^* + \underline{L}_f\overline{V})/\sigma = B\).
\Halmos
\endproof
\proof{Proof for Remark \ref{rmk1}.} Let \(x:=x_{\hat j}\), \(x^*:=x_{\hat j}^*\), and \(y:=y_{\hat j}^*\) be the quantities provided by Theorem \ref{thm8}. Since \(X=\mathbb{R}^n\), we have \(N_X(x^*)=\{0\}\). By Theorem \ref{thm8} and Lemma \ref{lem7}, \(\|x-x^*\|^2\le \bar\delta\), \(\|\nabla f(x^*)+\sum_{i=1}^m y_i\nabla h_i(x^*)\|^2\le \bar\epsilon\), \(\sum_{i=1}^m |y_i h_i(x^*)|\le \bar\epsilon\), and \(\|y\|_1\le B\). Lipschitz continuity of \(\nabla f\) and \(\nabla h_i\) gives \(\|\nabla f(x)+\sum_{i=1}^m y_i\nabla h_i(x)-\nabla f(x^*)-\sum_{i=1}^m y_i\nabla h_i(x^*)\| \le (\underline L_f + B\max_{i\in[m]}L_{h_i})\|x-x^*\|\), hence \(\|\nabla f(x)+\sum_{i=1}^m y_i\nabla h_i(x)\|^2=O(\bar\delta+\bar\epsilon)\). Likewise, Lipschitz continuity of \(h_i\) yields \(\sum_{i=1}^m |y_i h_i(x)| \le \sum_{i=1}^m |y_i h_i(x^*)| + \sum_{i=1}^m y_i|h_i(x)-h_i(x^*)| \le \bar\epsilon + B\max_{i\in[m]}M_{h_i}\|x-x^*\| = O(\sqrt{\bar\delta}+\bar\epsilon)\). Finally, since \(h_i(x^*)\le 0\) and each \(h_i\) is \(M_{h_i}\)-Lipschitz continuous, \([h_i(x)]_+ \le |h_i(x)-h_i(x^*)| \le M_{h_i}\|x-x^*\|\), so \(\max_{i\in[m]}[h_i(x)]_+=O(\sqrt{\bar\delta})\). \Halmos

\proof{Proof of Theorem \ref{thm8}.}
For each proximal subproblem, Algorithm~\ref{alg:inexact proximal} returns $x_j$ with objective suboptimality and constraint violation at most $\tilde\epsilon$. Lemma~\ref{lem6} then gives
\[
\frac{\underline L_f}{2}\|x_{\hat j}^*-x_{\hat j}\|^2
\le
f(x_{\hat j};x_{\hat j-1})-f(x_{\hat j}^*;x_{\hat j-1})
+\langle y_{\hat j}^*, h(x_{\hat j})\rangle
\le (1+B)\tilde\epsilon,
\]
so $\|x_{\hat j}^*-x_{\hat j}\|^2 \le \tfrac{2(1+B)\tilde\epsilon}{\underline L_f},$
and then $x_{\hat j}$ lies within distance $O(\sqrt{\tilde\epsilon/\underline L_f})$ of the exact KKT point $x_{\hat j}^*$ of the selected convex subproblem. Evaluating Lemma~\ref{lem6} at $x_{\hat j-1}$ and using the choice of $\hat j$ yields
\[
\left[d\left(\nabla f(x^*_{\hat j})+\sum_{i=1}^m y^*_{\hat j,i}\nabla h_i(x^*_{\hat j}),-N_X(x^*_{\hat j})\right)\right]^2
\le \frac{8\underline L_f M_fD_X}{J} + (1+B)\tilde\epsilon.
\]
Indeed, since \(f\) is \(M_f\)-Lipschitz continuous on the compact set \(X\), for every run of the algorithm we have
\[
f(x_0)-f(x_J)\le |f(x_0)-f(x_J)| \le M_f\|x_0-x_J\| \le M_fD_X.
\]
Therefore, if
\[
J \ge \frac{16\underline L_fM_fD_X}{\bar\epsilon},
\qquad
\tilde\epsilon \le \min\left\{\frac{\bar\epsilon}{2(1+B)},\frac{\underline L_f\bar\delta}{2(1+B)}\right\},
\]
then $x_{\hat j}^*$ is a $\bar\epsilon$-KKT point and $x_{\hat j}$ is an \((\bar\epsilon,\bar\delta)\)-near-KKT point of \eqref{cvxctr-model}. Each convex subproblem is solved by LCG in \(\mathcal O(\tilde\epsilon^{-2}\log(\tilde\epsilon^{-1}))\) CGO iterations, so multiplying by $J=O(\bar\epsilon^{-1})$ gives the stated total complexity. \Halmos
\endproof

\proof{Proof of Corollary \ref{cor:ipp-lcg-eps}.}
Set \(\bar\epsilon=\bar\delta=\epsilon\) in Theorem \ref{thm8}. Then \(\mathcal{O}\!\left(\bar{\epsilon}^{-1}\max\{\bar{\delta}^{-2}\log(\bar{\delta}^{-1}),\bar{\epsilon}^{-2}\log(\bar{\epsilon}^{-1})\}\right)=\mathcal{O}(\epsilon^{-1}\cdot \epsilon^{-2}\log(\epsilon^{-1}))=\mathcal{O}(\epsilon^{-3}\log(\epsilon^{-1}))\).
\Halmos
\endproof

\section{Auxiliary Numerical Studies}

\subsection{Portfolio Selection Models}\label{appen-num-portf}
\paragraph{Risk formulations.}
The objective in \eqref{portf-var-sa} is a discontinuous step-risk. For the smooth nonconvex model \eqref{portf-nonconvex}, we replace it by the sigmoid surrogate \(\tilde{f}_{\theta}(x)= K^{-1}\sum_{k=1}^K (1 + \exp\{ (-R_k + r_k^\top x)/\theta\})^{-1}\), which converges pointwise to the step-risk as \(\theta\to 0\). For the convex model \eqref{portf-convex}, write the benchmark-underperformance loss as \(L(x):=R-r^\top x\), so that the original objective is \(\mathbb{E}[\mathbb{I}\{L(x)>0\}]\). A standard convex surrogate is to replace this discontinuous probability by the tail-risk functional \(\CVaR_{\alpha}(L(x))\), which admits the Rockafellar--Uryasev variational representation \(\CVaR_{\alpha}(L(x))=\inf_{u\in\mathbb{R}}\{u+\alpha^{-1}\mathbb{E}[L(x)-u]_+\}\);
see \citet{rockafellar2000optimization,rockafellar2002conditional}. Replacing the expectation by the empirical average gives exactly the objective in \eqref{portf-convex}. This objective is convex in \((x,u)\) because each term \([R_k-r_k^\top x-u]_+\) is the hinge of an affine function, and the full feasible region remains convex after adding the simplex constraint, the box constraints on \((u,v)\), and the cap-type sparsity surrogate \eqref{portf-sparsity-surrogate}. In implementation, we further restrict the auxiliary variable \(u\) to a data-dependent interval \([\underline u,\bar u]\) to keep the feasible set compact.

\paragraph{Sparsity controls in the portfolio training models.}
The exact cardinality budget \(s(x)\le \Psi\) is used in the paper as the motivating sparse-portfolio requirement, but the projection-free training problems used in the portfolio experiments replace this combinatorial restriction by the convex cap-type surrogate \(g_{\Psi}(x,v):=Nv+\Psi^{-1}\sum_{i=1}^N [x_i-v]_+-N/\Psi\le 0\), together with \(\underline v\le v\le \bar v\).
Here \(v\) is an auxiliary threshold variable optimized jointly with \(x\). This does not reproduce the exact support constraint, but it discourages weight from being spread across many assets while remaining convex in \((x,v)\), since each hinge term is convex and the remaining terms are affine. The functional \(g_{\Psi}\) is convex in \((x,v)\) because each hinge term \([x_i-v]_+\) is convex. In the portfolio experiments, we calibrate the surrogate by setting \(\underline v=10^{-4}\) and \(\bar v=1/\Psi\); the same \(10^{-4}\) threshold is also used when reporting realized support, and the offset \(N/\Psi\) makes the equal-weight \(\Psi\)-asset portfolio satisfy \(g_{\Psi}(x,v)=0\) at \(v=1/\Psi\). This keeps the portfolio training formulations within the functional-constrained template studied in the main text. Because our numerical comparison is ultimately concerned with portfolio quality and realized sparsity, the main-text tables report the actual support \(s(x)\) of the selected portfolio rather than the value of \(g_{\Psi}(x,v)\).


\paragraph{Sparsity targets and evaluation.}
In the main-text portfolio experiments (Section~\ref{sec:numerical-portf}), we evaluate each method as a sparse portfolio solver under a target support size \(\Psi\). For the projection-free methods, \(\Psi\) is also the calibration parameter used in the convex surrogate sparsity constraint \(g_{\Psi}(x,v)\le 0\) of the training models. We report the best recorded iterate that reaches \(s(x)\le \Psi\), and we penalize methods that never reach the target.

\paragraph{Main-text evaluation protocol.}\label{append:portf-eval}
All methods are run for the same 5-second wall-clock budget, and for each dataset every method uses the same 70/30 train/test split. Among all recorded iterates with \(s(x)\le \Psi\), we select the one with minimum training \(\mathrm{StepRisk}_{\delta}\), breaking ties by smaller training \(\CVaR_{\alpha}\), then by larger support (closer to \(\Psi\)), then by earlier time. Here \(\CVaR_{\alpha}(x):=\min_{u\in\mathbb{R}}\{u+(\alpha K_{\mathrm{train}})^{-1}\sum_{k\in \mathrm{train}}[R_k-r_k^\top x-u]_+\}\), with \(\alpha=0.1\). If a method never reaches \(s(x)\le \Psi\), we instead select the iterate minimizing the penalized score \(\mathrm{StepRisk}_{\delta} + \rho\, [s(x)-\Psi]_+/\Psi\) with \(\rho=1\).
%
{\color[rgb]{0,0,0.8}
\begin{table}[tp]
\centering
\caption{Per-dataset iteration counts under the 5-second portfolio time-budget rerun. For LCG-style and projected first-order methods, the entries are direct iteration counts. For Greedy+Refit, the entries are refit iterations after one greedy selection pass. For IRL1 and PenPGD, the entries are total inner updates accumulated across all candidate subruns within the same 5-second budget.}
\resizebox{0.75\columnwidth}{!}{%
\begin{tabular}{l|ccccccccc}
\toprule
 & LCG & CoexDurCG & FISTA & Greedy+Refit & IRL1 & PenPGD & PGD & PGD-IHT & MD(ent) \\
\midrule
DJ & 6989 & 41893 & 53015 & 34319 & 34931 & 32519 & 33568 & 29030 & 39399 \\
FF49 & 4838 & 28146 & 36190 & 28762 & 23295 & 23603 & 22769 & 24210 & 28755 \\
ND100 & 11910 & 44150 & 52479 & 37459 & 36152 & 34183 & 34682 & 31601 & 41650 \\
FTSE100 & 11029 & 44047 & 44434 & 36578 & 34471 & 34638 & 34620 & 30273 & 40716 \\
SP500 & 8033 & 23857 & 27688 & 37269 & 24316 & 23536 & 24249 & 21514 & 29636 \\
NDComp & 14 & 23 & 20 & 31473 & 50 & 46 & 14 & 15 & 15 \\
\bottomrule
\end{tabular}%
}
\vspace{0.5ex}
\resizebox{0.6\columnwidth}{!}{%
\begin{tabular}{l|ccccccc}
\toprule
 & IPP-LCG & Greedy+Refit & IRL1 & PenPGD & PGD & PGD-IHT & MD(ent) \\
\midrule
DJ & 24356 & 73508 & 40537 & 39549 & 38563 & 35295 & 45516 \\
FF49 & 17706 & 56327 & 20516 & 19630 & 19073 & 17932 & 19984 \\
ND100 & 25819 & 74408 & 43211 & 43118 & 42960 & 37783 & 51501 \\
FTSE100 & 24909 & 73337 & 36220 & 35501 & 34480 & 32208 & 43249 \\
SP500 & 15463 & 72400 & 16082 & 15843 & 15091 & 14134 & 16688 \\
NDComp & 11 & 45499 & 44 & 51 & 18 & 16 & 17 \\
\bottomrule
\end{tabular}%
}
\label{tab:portf-time-iters}
\end{table}
}
\subsection{Additional Portfolio Baseline Implementations}\label{sec:append:portf-baselines}
This subsection summarizes the non-projection-free portfolio baselines used in the main-text tables and how we adapt them to the common time-budgeted sparse-solver comparison in Section~\ref{sec:numerical-portf}. The projection-free comparator CoexDurCG is the method of \citet{lan2020conditional}; the remaining baselines are standard projected first-order methods or sparse-optimization heuristics instantiated for the portfolio training models in Section~\ref{sec:app-portf}.
\begin{itemize}
    \item \textbf{PGD / MD.} Projected gradient descent is the basic Euclidean first-order baseline for convex-constrained optimization \citep{bertsekas2015convex}, while mirror descent with the entropy proximity is its natural non-Euclidean counterpart on the simplex \citep{nemirovskij1983problem,nemirovski2004prox}. We run both directly on the portfolio training objectives within the same 5-second budget. For the convex \CVaR\ model, the auxiliary variable \(u\) is updated by projected gradient on \([\underline u,\bar u]\); for the sparsity-threshold variable \(v\), we analogously use projected Euclidean updates on \([\underline v,\bar v]\).
    \item \textbf{PGD-IHT and FISTA.} PGD-IHT combines projected-gradient steps with iterative hard-thresholding-style support truncation \citep{blumensath2008iterative}: after each gradient step, we keep the largest \(\Psi\) components of \(x\) and renormalize on the simplex. FISTA \citep{beck2009fast} is used as an accelerated projected first-order baseline for the convex model only; specifically, we apply it to a softplus-smoothed version of the \CVaR\ objective over \((x,u)\in X\times[\underline u,\bar u]\), since the original hinge objective is nonsmooth.
    \item \textbf{Greedy+Refit and IRL1.} Greedy+Refit is a one-pass support-selection heuristic inspired by greedy sparse approximation methods such as orthogonal matching pursuit \citep{tropp2007signal}: one gradient/subgradient evaluation is used to rank assets and choose a candidate support, after which we refit the portfolio on that support. IRL1 uses iteratively reweighted \(\ell_1\) surrogates to promote sparsity \citep{candes2008enhancing}; in our implementation, a small grid of regularization parameters is explored within the same 5-second budget, and the best iterate is selected using the rule in Appendix~\ref{append:portf-eval}.
    \item \textbf{PenPGD.} PenPGD applies projected-gradient updates to penalized training objectives with several standard nonconvex sparsity penalties: reweighted-\(\ell_1\)/log-sum-type surrogates \citep{candes2008enhancing}, {\sf SCAD} \citep{fan2001variable}, {\sf MCP} \citep{zhang2010nearly}, and capped-\(\ell_1\)/multi-stage convex relaxation surrogates \citep{zhang2010analysis}. We run these candidate penalties under the same 5-second budget and report the best iterate under the common selection rule. As noted in Table~\ref{tab:portf-time-iters}, the reported update counts for IRL1 and PenPGD aggregate the inner updates performed across these candidate subruns.
\end{itemize}

\subsection{Portfolio Dataset Statistics and Reference Risks}\label{sec:append:portf-results}
This subsection collects the supporting dataset statistics and reference risk values used in the portfolio experiments (Tables~\ref{tab:stock-data}--\ref{tab:portf-ref-risk}).
\begin{table}[tp]
\caption{Features of the stock market dataset.}
\centering
{\small
\begin{tabular}{llccc}
\toprule
Instance & Description & $\#$ of assets ($N$) & $\#$ of weeks ($K$) & Cardinality ($\Psi$)\\
\midrule
DJ & Dow Jones Industrial Average (USA) & 28 & 1363 & 5 \\
FF49 & Fama and French 49 Industry (USA) & 49 & 2325 & 9\\
ND100 & NASDAQ 100 (USA) & 82 & 596 & 16 \\
FTSE100 & FTSE 100 (UK) & 83 & 717 & 16 \\
SP500 & S\&P 500 (USA) & 442 & 595 & 22 \\
NDComp & NASDAQ Composite (USA) & 1203 & 685 & 60\\
\bottomrule
\end{tabular}}
\label{tab:stock-data}
\end{table}

\begin{table}[tp]
\caption{Reference out-of-sample step-risk values at threshold $\delta=0.0005$ for two simple long-only portfolios under the same 70/30 time split: all-cash ($x=0$; infeasible under $\sum_i x_i = 1$ but included as a reference) and equal-weight ($x_i=1/N$).}
\centering
{\small
\begin{tabular}{lcc}
\toprule
Instance & StepRisk$_{\delta}$(test), cash & StepRisk$_{\delta}$(test), equal-weight\\
\midrule
DJ & 0.5232 & 0.3472\\
FF49 & 0.6003 & 0.0029\\
ND100 & 0.5754 & 0.4134\\
FTSE100 & 0.5231 & 0.2778\\
SP500 & 0.5922 & 0.3855\\
NDComp & 0.5728 & 0.4369\\
\bottomrule
\end{tabular}}
\label{tab:portf-ref-risk}
\end{table}

%
%

\subsection{IMRT Models}\label{sec:appen-imrt-cvx}
For each patient, the anatomy is discretized into voxels $v\in \mathcal V$. Each beam angle $a\in A$ has a large implicit set of deliverable apertures $E_a$, and $y_{a,e}$ denotes the intensity assigned to aperture $e\in E_a$. If $x_{lr}^{a,e}$ indicates whether beamlet $(l,r)$ is open in aperture $e$, then the dose delivered to voxel $v$ is \(z_v = \sum_{a\in A}\sum_{e\in E_a} \sum_{l=1}^m\sum_{r= 1}^n RD_{(l,r) v}x^{a,e}_{lr}y_{a,e}\).
We use $K_u$ and $K_o$ to index underdose and overdose criteria, with $S_k$ the voxel set in criterion $k$, $N_k=|S_k|$, and $(b_k,p_k)$ the corresponding clinical threshold and tolerated fraction.

\paragraph{Convex Formulation.}
Following \citet{lan2020conditional}, we replace the dose-volume constraints by \CVaR-type upper bounds and solve
\begin{align}
	\min\ f(z):= \frac{1}{N_v}&\sum\limits_{v\in \mathcal{V}}\underline{w}_v \left[\underline{T}_v - z_v \right]^{2}_{+} + \bar{w}_v\left[ z_v - \bar{T}_v\right]^{2}_{+} \label{cvx-obj}\\
	&\st \ -\tau_k + \frac{1}{p_kN_k}\sum\limits_{v\in S_k} \left[\tau_k - z_v \right]_{+} \le -b_k,\ \forall k \in K_u, \label{cvx-ctr1} \\
	& \tau_k + \frac{1}{p_kN_k}\sum\limits_{v\in S_k} \left[z_v - \tau_k \right]_{+} \le b_k,\ \forall k \in K_o, \label{cvx-ctr2} \\
	& z_v = \sum\limits_{a\in A}\sum\limits_{e\in E_a} \sum\limits_{l=1}^m\sum\limits_{r= 1}^n RD_{(l,r) v}x^{a,e}_{lr}y_{a,e},\ \forall v\in S_k,\ k\in K_u\cup K_o,\\
	& \sum\limits_{a\in A} \max_{e\in E_a} y_{a,e} \le \Phi,  \label{cvx-ctr3} \\
	& \sum\limits_{a\in A}\sum\limits_{e\in E_a} y_{a,e} \le 1,  \label{cvx-ctr4} \\
	& y_{a,e} \ge 0,\ \forall a\in A,\ e\in E_a,\\
	& \underline{\tau}_k \le \tau_k \le \bar{\tau}_k,\qquad k\in K_u \cup K_o. \label{cvx-last-ctr}
\end{align}
Here \eqref{cvx-ctr3} is the group-sparsity control, and the nonsmooth constraints \eqref{cvx-ctr1}--\eqref{cvx-ctr3} fit the structured form \eqref{eq:ns-structure-h}, so they are handled by the smoothed CGO oracle in Algorithm~\ref{alg:CG oracle nonsmooth}.

\paragraph{Nonconvex Formulation.}
The nonconvex model keeps the original dose-volume criteria in the objective and retains only the deliverability and group-sparsity constraints:
\begin{equation}
\begin{gathered}
	\min\  f(z) := \sum\limits_{k\in K_u} \frac{w_k}{N_k}\sum\limits_{v\in S_k} \mathbb{I}_{\{z_v < \tau_k\}} + \sum\limits_{k\in K_o} \frac{w_k}{N_k}\sum\limits_{v\in S_k} \mathbb{I}_{\{z_v > \tau_k\}}\\
	\st\  z_v = \sum\limits_{a\in A}\sum\limits_{e\in E_a} \sum\limits_{l=1}^m\sum\limits_{r= 1}^n RD_{(l,r) v}x^{a,e}_{lr}y_{a,e},\ \forall v\in S_k,\ k\in K_u\cup K_o,\\
	 \sum\limits_{a\in A} \max_{e\in E_a} y_{a,e} \le \Phi, \\
	 \sum\limits_{a\in A}\sum\limits_{e\in E_a} y_{a,e} \le 1, \\
	 y_{a,e} \ge 0.
\end{gathered}
\end{equation}
We solve this model through the sigmoid approximation described in the main text:
\begin{equation*}
\begin{aligned}
\tilde{f}^k_{\theta}(x)
&= \frac{1}{N_k}\sum\limits_{v\in S_k}\frac{1}{1 + \exp\{ (z_v - \tau_k)/\theta\}},
&& k\in K_u,\\
\tilde{f}^k_{\theta}(x)
&= \frac{1}{N_k}\sum\limits_{v\in S_k}\frac{1}{1 + \exp\{ (-z_v + \tau_k)/\theta\}},
&& k\in K_o.
\end{aligned}
\end{equation*}
As $\theta \to 0$, these surrogates converge pointwise to the corresponding indicator losses.

\subsection{IMRT Synthetic Dataset}\label{sec:imrt-synthetic}
The synthetic instances mimic the scale and geometry of a real IMRT case. Each instance has 180 beam angles, associated dose matrices, and three clinical criteria (two underdose and one overdose). Instances 1--2 use a coarser beamlet discretization, while Instances 3--4 use a finer discretization and are therefore substantially larger.
\begin{table}[h]
\caption{Features of the synthetic dataset}
	\centering
	{\small
	\begin{tabular}{ccccccc}
		\toprule
		 Instance & $\#$ of angles &	$\#$ of voxels & $\#$ of beamlets  & Accuracy &$b_k$ & $p_k$ \\
		\midrule
	1	& 180 &$4096$ &  100 & 1.0  & $[40, 50, 100]$  &$[0.01, 0.01, 0.05]$\\
	2	 & 180 &$4096$&  100 & 1.0 & $[50, 60, 80]$  &$[0.01, 0.01, 0.01]$\\
	3	& 180 &$262144$ &  2000 & 0.25 & $[40, 50, 100]$  &$[0.01, 0.01, 0.05]$\\
	4	& 180 &$262144$ &  2000 & 0.25 &$[50, 60, 80]$  &$[0.01, 0.01, 0.01]$\\
		\midrule
		\bottomrule
	\end{tabular}}
	\label{tab:ran-data}
\end{table}

\end{document}